\RequirePackage{amsmath}
\documentclass[runningheads]{llncs}
\usepackage[T1]{fontenc}
\usepackage{newtxtext}       %
\usepackage[varvw]{newtxmath}       

\usepackage{graphicx}
\usepackage{mathtools}
\usepackage{dsfont}
\usepackage{enumitem}
\usepackage{xcolor}
\usepackage{tabularx}
\usepackage{booktabs}
\usepackage{tikz}
\usetikzlibrary{arrows.meta}
\colorlet{wolfecol}{blue!70!black}
\definecolor{grayif}{gray}{0.6}
\usepackage{pdflscape}
\usepackage{afterpage}

\usepackage[colorinlistoftodos]{todonotes}
\expandafter\let\csname algorithm*\endcsname\relax
\expandafter\let\csname endalgorithm*\endcsname\relax

\usepackage[ruled,vlined]{algorithm2e}
\usepackage{setspace}
\newcommand\mycommfont[1]{\footnotesize\ttfamily\textcolor{gray}{#1}}
\SetCommentSty{mycommfont}
\newcommand{\vgap}{\vspace{.2em}}
\SetKw{KwAnd}{and}
\SetKw{KwOr}{or}
\SetKw{KwNot}{not}
\DontPrintSemicolon
\SetNlSty{}{\color{black}\scriptsize}{}
\SetArgSty{textup}
\newcommand{\algstrut}{\vphantom{$\ivec{d}{B\setminus\{j\},j}\ivec{\bar{\lambda}}{S_0}\ivec{x}{B_0}\vec{w_S}\R^D\Q^D\smash[b]{\rc^B_j}$}}
\DeclareRobustCommand{\algcapstrut}{\vphantom{$(A,\vec{b},\vec{c},B_0,\ivec{p}{n},S_0)$}}
\SetKw{KwReturn}{return}
\SetKwFor{For}{for}{\!:}{endfor}
\SetKwFor{While}{while}{\!:}{endw}
\SetKwFor{ForEach}{for each}{\!:}{endfor}
\SetKwIF{If}{ElseIf}{Else}{if}{\!:}{else if}{else}{endif}
\SetKwComment{tcp}{(}{)}

\iftrue
\makeatletter
\renewcommand{\SetKwInOut}[2]{%
  \sbox\algocf@inoutbox{\KwSty{#2}\algocf@typo:}%
  \expandafter\ifx\csname InOutSizeDefined\endcsname\relax
    \newcommand\InOutSizeDefined{}\setlength{\inoutsize}{\wd\algocf@inoutbox}%
    \sbox\algocf@inoutbox{\parbox[t]{\inoutsize}{\KwSty{#2}\algocf@typo:\hfill}~}\setlength{\inoutindent}{\wd\algocf@inoutbox}%
  \else
    \ifdim\wd\algocf@inoutbox>\inoutsize%
    \setlength{\inoutsize}{\wd\algocf@inoutbox}%
    \sbox\algocf@inoutbox{\parbox[t]{\inoutsize}{\KwSty{#2}\algocf@typo:\hfill}~}\setlength{\inoutindent}{\wd\algocf@inoutbox}%
    \fi%
  \fi
  \algocf@newcommand{#1}[1]{%
    \ifthenelse{\boolean{algocf@inoutnumbered}}{\relax}{\everypar={\relax}}%
    {\let\\\algocf@newinout\hangindent=\inoutindent\hangafter=1\parbox[t]{\inoutsize}{\KwSty{#2}\algocf@typo:\hfill}~##1\par}%
    \algocf@linesnumbered
  }}%
\makeatother
\fi
\SetKwInOut{Input}{input}
\SetKwInOut{Output}{output}

\usepackage{cite}
\usepackage[hypertexnames=false,hidelinks]{hyperref}
\usepackage{orcidlink}

\renewcommand{\vec}[1]{\boldsymbol{#1}}
\newcommand{\N}{\mathds{N}}
\newcommand{\Z}{\mathds{Z}}
\newcommand{\Q}{\mathds{Q}}
\newcommand{\R}{\mathds{R}}
\newcommand{\one}{\mathbf 1}
\newcommand{\zero}{\mathbf 0}
\newcommand{\norm}[1]{\lVert#1\rVert}
\newcommand{\rc}{\bar{c}}
\newcommand{\rg}{\bar{\gamma}}
\DeclareMathOperator*{\polytope}{\mathcal P}

\DeclareMathOperator*{\argmin}{arg\,min}

\DeclareMathOperator{\conv}{conv}
\DeclareMathOperator{\aff}{aff}

\newcommand{\ivec}[2]{\vec{#1}^{\scriptscriptstyle(#2)}}
\newcommand{\icoord}[3]{#1^{\scriptscriptstyle(#2)}_{#3}}

\newif\ifsubmission
\submissionfalse

\newcommand{\proofinappendix}{\leavevmode\hskip-\labelsep\footref{fn:proofs}\ }
\newcommand{\corollaryinappendix}{\leavevmode\hskip-\labelsep\footref{fn:corollaries}\ }
\newcommand{\tabeq}[1]{\refstepcounter{equation}(\theequation)\label{#1}}
\spnewtheorem{claimn}{Claim}{\itshape}{\rmfamily}
\spnewtheorem{openproblem}{Open Problem}{\bfseries}{\itshape}

\begin{document}
\title{Wolfe's Minimum-Norm-Point Method\\ Generalizes the Simplex Method}
%
%
\ifsubmission
\author{Anonymous Authors}
\authorrunning{Anonymous Authors}
\institute{}
\else
\author{Yann Disser\,\orcidlink{0000-0002-2085-0454} \and Nils Mosis\,\orcidlink{0000-0002-0692-0647}}
\authorrunning{Y. Disser and N. Mosis}
%
\institute{TU Darmstadt, Germany\\
\email{\{disser,mosis\}@mathematik.tu-darmstadt.de}}

\fi
\maketitle              
\begin{abstract}
Wolfe's method computes the point of minimum Euclidean norm in the convex hull of finitely many points.
The method is parameterized by an \emph{insertion rule} and it has been an open problem whether Wolfe's original rule guarantees a polynomial number of iterations (De Loera, Haddock, and Rademacher; STOC'18).
We resolve this question by proving an exponential lower bound.
In fact, we establish a tight connection to the simplex method for linear programming that allows us to transfer many lower bounds for simplex pivot rules to corresponding insertion rules of Wolfe's method.
These include not only Wolfe's linopt rule via Dantzig's pivot rule, but also many other classical, randomized and history-based rules.
As a corollary, we obtain the first exponential example for the Lawson–Hanson algorithm for nonnegative least squares, as well as PSPACE-hardness of deciding whether the linopt rule ever uses a given point.


\end{abstract}

\section[Introduction]{Introduction}
\label{sec:intro}

In 1974, Philip Wolfe~\cite{Wolfe74,Wolfe76} introduced a combinatorial active-set method for the \emph{minimum-norm-point problem}
\begin{equation}
    \min\ \norm{\vec{w}},\quad\text{subject to}\quad\vec{w}\in\conv\{\ivec{p}{1},\dots,\ivec{p}{n}\},
    \label{eq:mnp}
    \tag{MNP}
\end{equation}
which asks for the point of smallest Euclidean norm in the convex hull of finitely many points~$\ivec{p}{1},\dots,\ivec{p}{n}\in\Q^D$.
Wolfe's method maintains an iterate~$\vec{w}$ as well as an active set~$\mathcal{S} \subseteq \{\ivec{p}{1},\dots,\ivec{p}{n}\}$ of at most~$D+1$ affinely independent points.

The method repeatedly computes the minimum-norm point~$\vec{y}$ in the affine hull~$\aff(\mathcal{S})$ and moves~$\vec{w}$ in a straight line towards~$\vec{y}$.
If~$\vec{w}$ hits a face of~$\conv{\mathcal{S}}$, it eliminates a point in~$\mathcal{S}$ outside this face and repeats the process.
Otherwise, it adds any point~$\ivec{p}{i}$ to~$\mathcal{S}$ that allows further progress towards the origin, i.e., with~$\vec{w}^\top \ivec{p}{i} < \norm{\vec{w}}^2$ and thus, for sufficiently small~$\varepsilon > 0$, 
$$\norm{(1-\varepsilon)\vec{w} + \varepsilon\ivec{p}{i}}^2 = \norm{\vec{w}}^2
+ 2\varepsilon\bigl(\vec{w}^\top\ivec{p}{i}
                    - \norm{\vec{w}}^2\bigr)
+ \varepsilon^2\norm{\ivec{p}{i}-\vec{w}}^2 < \norm{\vec{w}}^2.$$
The method terminates with the global solution~$\vec{w}$ once it does not find any point to add to~$\mathcal{S}$.
See Figure~\ref{fig:wolfe-example} for an illustrated example.

\begin{figure}[tbp]
    \centering
    \begingroup
    \definecolor{edgecol}{gray}{0.55}
    \scalebox{0.9}{%
    \begin{tikzpicture}[x=.75cm,y=.75cm,
        trajectory/.style={wolfecol,line width=1.2pt,-{Stealth[length=2.3mm,width=1.7mm]},shorten <=2.5pt,shorten >=3pt},
        every node/.style={font=\small}]
      \coordinate (O) at (0,0);
      \coordinate (p1) at (-7,5.25);
      \coordinate (p2) at (5,4.25);
      \coordinate (p3) at (-5,1.25);
      \coordinate (p4) at (1,1.25);
      \coordinate (w1) at ({56/145},{672/145});
      \coordinate (z1) at ({55/234},{110/39});
      \coordinate (w2) at ({-165/218},{275/109});
      \coordinate (z2) at (-.375,1.25);
      \coordinate (w3) at (0,1.25);

      \draw[black!20,-{Stealth[length=1.5mm]},line width=.5pt] (-7.6,0)--(5.7,0);
      \draw[black!20,-{Stealth[length=1.5mm]},line width=.5pt] (0,-.22)--(0,5.9);

      \fill[black!3] (p1)--(p2)--(p4)--(p3)--cycle;
      \draw[edgecol,line width=.7pt] (p1)--(p2)--(p4)--(p3)--cycle;
      \draw[edgecol,line width=.6pt,dash pattern=on 2pt off 2.5pt] (p2)--(p3);

      \foreach \a/\b in {p1/w1,w1/z1,z1/w2,w2/z2}
        \draw[trajectory] (\a)--(\b);
      \draw[trajectory,-{Stealth[length=.9mm,width=.75mm]},shorten <=2.1pt,shorten >=2.7pt] (z2)--(w3);

      \foreach \p in {p1,p2,p3,p4}\fill (\p) circle[radius=2.2pt];
      \foreach \p in {p1,w1,w2}\fill[wolfecol] (\p) circle[radius=2.5pt];
      \foreach \p in {z1,z2}\draw[wolfecol,line width=.9pt,fill=white] (\p) circle[radius=2pt];
      \draw[wolfecol,line width=.8pt,fill=white] (w3) circle[radius=2.6pt];
      \fill[wolfecol] (w3) circle[radius=1.4pt];
      \fill (O) circle[radius=1.7pt];

      \node[anchor=south east,xshift=1pt,yshift=1pt,inner sep=0pt] at (p1) {$\ivec{p}{1}$};
      \node[above right=3pt,inner sep=1pt] at (p2) {$\ivec{p}{2}$};
      \node[below left=1pt,inner sep=1pt] at (p3) {$\ivec{p}{3}$};
      \node[below right=1pt,inner sep=1pt] at (p4) {$\ivec{p}{4}$};
      \node[below left=1pt,inner sep=1pt] at (O) {$\vec{0}$};

    \end{tikzpicture}}
    \endgroup
    \caption{Wolfe's method with the linopt rule, starting at~$\ivec{p}{1}$. 
    Filled circles mark the positions at which the points~$\ivec{p}{2}$,~$\ivec{p}{3}$, and~$\ivec{p}{4}$ are inserted into the active set, in this order.
    At the open circles,~$\ivec{p}{1}$ and~$\ivec{p}{2}$ are deleted, respectively.}
    \label{fig:wolfe-example}
\end{figure}

The behavior of Wolfe's method crucially depends on the \emph{insertion rule} that governs the choice which point to add to the active set~$\mathcal{S}$ if there is more than one candidate.
It has been an open question whether an insertion rule exists that guarantees a polynomial number of iterations.
De~Loera, Haddock, and Rademacher~\cite{DHR} showed that Wolfe's method with the \emph{minnorm} rule that chooses a point~$\ivec{p}{i}$ minimizing~$\norm{\ivec{p}{i}}$ to enter~$\mathcal{S}$ can take an exponential number of iterations.
They asked for such examples for other insertion rules, in particular for Wolfe's original \emph{linopt} rule minimizing~$\vec{w}^\top \ivec{p}{i} - \norm{\vec{w}}^2$, but their question has remained open~\cite[Section~1.1]{FKV25}.
Notably, the linopt rule has a significant advantage over the minnorm rule: It can be operated if the convex hull of input points is given in exterior representation (i.e., as inequalities) rather than in interior representation (i.e., by its vertices), see the discussion in~\cite{DHR}.

The existence of efficient insertion rules has been an open problem for active-set methods beyond Wolfe's method.
Recently, expanding on prior work by Disser and Mosis~\cite{DM25}, Bach et al.~\cite{BDHM26} showed that no such rule exists for the classical active-set method with a convex quadratic \emph{maximization} (or equivalently, \emph{concave} quadratic minimization) objective, explicitly leaving the case of convex quadratic minimization open.
Wolfe’s method is an active-set method for convex quadratic minimization, since the minimum-norm-point problem (MNP) can be written as the convex quadratic problem
\begin{align}
\min_{\lambda\in\R^n}\ \frac12\sum_{i,j=1}^n\lambda_i\lambda_j\,(\ivec{p}{i})^\top\ivec{p}{j}
\quad\text{subject to}\quad \sum_{i=1}^n\lambda_i=1,\ \lambda_i\geq0\ \text{for all }i.\label{eq:CQM}
\end{align}

Even more importantly, linear programming reduces to the minimum-norm-point problem~\cite{FHI06}, even in strongly polynomial time~\cite{DHR}.
Wolfe's algorithm being a combinatorial active-set method means that an efficient insertion rule would likely yield a strongly polynomial algorithm for linear programming, which would resolve Smale's ninth problem~\cite{Sma98} and constitute a major breakthrough in the field of discrete optimization.
Since its inception by Dantzig in 1947 \cite{Dan82}, the simplex method has been the most prominent candidate for an efficient combinatorial algorithm.
Much like Wolfe's method, it is governed by a selection (or \emph{pivot}) rule, and exponential or superpolynomial lower bounds have been found for a multitude of rules~\cite{KM72,Jer73,AC78,GS79,Mur80,Gol83,AZ98,Fri11,FHZ11,FHZ14,AF17,DM23,DFH23}
(see Table~\ref{tab:wolfe-rules} for some of those).
A polynomial pivot rule for the simplex method would not only answer Smale's problem but also prove the polynomial Hirsch conjecture \cite{Dan63,San12}.

In this paper, we establish a tight connection between Wolfe's method and the simplex method that allows us to transfer many results from the latter to the former, including an exponential lower bound for the linopt insertion rule.

\subsubsection*{Our results.}

Our main result is a reduction from linear programs (LPs) to minimum-norm-point problems (MNPs) of similar size with a one-to-one correspondence between executions of the simplex method and Wolfe's method.
In particular, every pivot rule for the simplex method turns into an insertion rule of Wolfe's method, and vice-versa.
The starting point of our result is a formulation, inspired by~\cite{DM25}, that makes the similarities of the methods apparent (Algorithms~\ref{alg:simplex} and~\ref{alg:wolfe}).
Formally, we prove the following.\footnote{\label{fn:lean}Theorem 1, Corollaries 1–4, and all numbered supporting statements have been formally verified in Lean 4: 
\url{https://doi.org/10.5281/zenodo.24136107}.}

\begin{theorem}\label{thm:intro-equivalence}
  Every bounded, non-degenerate LP with~$n$ variables and~$m$ equality constraints induces, in polynomial time, an MNP on~$n+1$ points in~$\Z^{m+1}$, such that the sequences of simplex pivots are in bijection with the sequences of Wolfe insertions.
\end{theorem}

Theorem~\ref{thm:intro-equivalence} establishes an intimate connection between the methods that has various consequences.
First, we obtain that Wolfe's method generalizes the simplex method in the sense that the existence of an insertion rule for Wolfe's method that turns it into a (strongly) polynomial algorithm implies the existence of a pivot rule with the same property for the simplex method.

\begin{corollary}\leavevmode\hskip-\labelsep\footnote{\label{fn:corollaries} We derive Corollaries~\ref{cor:intro-polynomial}--\ref{cor:intro-pspace} from Theorem~\ref{thm:equivalence} in Appendix~\ref{app:consequences}.}\ %
\label{cor:intro-polynomial}
    Every (strongly) polynomial insertion rule for Wolfe's method induces a (strongly) polynomial pivot rule for the simplex method.
\end{corollary}

Second, by understanding the translation between pivot and insertion rules in Theorem~\ref{thm:intro-equivalence}, we are able to transfer various lower bound examples for the simplex method to Wolfe's method.
In particular, every simplex pivot rule that bases its decision only on the reduced costs (up to a common factor), the indices, and the history of the run becomes a natural insertion rule:
Dantzig's rule\cite{Dan63} becomes linopt, Bland's rule\cite{Bla77} becomes least-index, and the rules of Zadeh~\cite{Zad80} and Cunningham~\cite{Cun79}, Random-Edge, and Random-Bland become natural memory-based and randomized insertion rules (Table~\ref{tab:wolfe-rules}).
We thus obtain lower bounds for each of these rules, in particular we close the open problem~\cite{DHR,FKV25} of the complexity for Wolfe's original linopt rule and recover the lower bound for the minnorm rule~\cite{DHR} (Appendix~\ref{app:cor-rules}).

\begin{corollary}\corollaryinappendix\label{cor:intro-rules}
    Wolfe’s method requires exponentially many iterations in the worst case with the linopt, minnorm, least-index, least-inserted, and round-robin rules, and superpolynomially many in expectation with the random and random-order rules.
    
    The lower bounds for the linopt and minnorm rules hold even when starting Wolfe's method at the input point of minimum norm.
\end{corollary}

Applied to the known reduction of least-distance programming to nonnegative least squares (NNLS)~\cite[Chapter~23]{LH74}, the Lawson--Hanson algorithm simulates Wolfe's method with the linopt rule~\cite[Section~1.1]{DHR} (see Appendix~\ref{app:nnls}), so Corollary~\ref{cor:intro-rules} also yields the first, to the best of our knowledge, exponential lower bound for this algorithm.

\begin{corollary}\corollaryinappendix\label{cor:intro-nnls}
    For every~$n\geq1$, there is a matrix~$E \in \Z^{(n+3)\times(2n+2)}$ such that the Lawson--Hanson algorithm for~$\min\{\norm{E\vec{\nu}-\vec{e_{n+3}}}:\vec{\nu}\geq\zero\}$ performs~$2^n+n+1$ insertions.
\end{corollary}

Finally, the construction behind Theorem~\ref{thm:intro-equivalence} allows to transfer the mightiness of the simplex method~\cite{FS15,DS18} to Wolfe's method.

\begin{corollary}\corollaryinappendix\label{cor:intro-pspace}
It is PSPACE-hard to decide whether Wolfe's method with the linopt rule ever inserts a given point, and it is NP-hard to determine the total number of insertions.
\end{corollary}

\subsubsection*{Related work.}
\label{sec:related}

The minimum-norm-point problem is a fundamental optimization problem that appears as a subproblem in many areas.
It computes projections onto polytopes and distances between polytopes, for example, in collision detection~\cite{GJK88} and in the training of support vector machines~\cite{KSBM00}.
If the origin is not in the convex hull of the points, then the minimum-norm point~$\vec{w^*}$ gives the solution~$\vec{w^*}/\norm{\vec{w^*}}^2$ of the least-distance problem $\min\{\norm{\vec{x}}:\vec{x}^\top\ivec{p}{j}\geq1\text{ for all }j\in[n]\}$, which Lawson and Hanson reduce to \emph{nonnegative least squares} (NNLS)~\cite[Chapter~23]{LH74}.
Applied to this reduction, the active-set algorithm of Lawson and Hanson for NNLS~\cite[Algorithm~(23.10)]{LH74}, which is behind MATLAB's \texttt{lsqnonneg} and SciPy's \texttt{nnls}, performs essentially the same steps as Wolfe's method with the linopt rule~\cite{DHR}.
With this rule, Wolfe's method is also the core of the Fujishige--Wolfe algorithm for submodular function minimization~\cite{Fuj80,FI11,FHI06}, and of the LP-Newton method for linear programming~\cite{FHI06,FHYZ09}.

Quantitative convergence bounds for Wolfe’s method with the linopt rule were established in~\cite{CJK14,LJJ15}. These bound the number of iterations needed to reach a prescribed accuracy in terms of numerical or geometric parameters of the input. For the closely related Lawson–Hanson algorithm, an exact termination bound is known that is polynomial in the dimension and a circuit imbalance measure~\cite{FKV25}. The dependence on these parameters can be exponential in the input encoding length, so these results do not imply polynomial worst-case running time.

\paragraph{Notation.}
\label{sec:notation}
We write~$[n]=\{1,\dots,n\}$.
For~$S\subseteq[d]$, the restriction of a vector~$\vec{x}\in\R^d$ to the coordinates in~$S$ is~$\vec{x_S}$, and a vector in~$\R^S$ is identified with its extension by zeros.
For a matrix~$A$, we denote column~$j$ by~$\vec{A_j}$ and the submatrix of the columns in~$S$ by~$A_S$.

\section{Preliminaries}
\label{sec:prelim}

\subsection{The simplex method}
\label{sec:simplex}

The simplex method solves linear programs of the form
\begin{equation}
    \min\ \vec{c}^\top\vec{x},\quad\text{subject to}\quad A\vec{x}=\vec{b},\ \vec{x}\geq\zero,
    \label{eq:lp}
    \tag{LP}
\end{equation}
with~$\vec{b}\in\Q^m$,~$\vec{c}\in\Q^n$, and~$A\in\Q^{m\times n}$ of rank~$m\geq1$.

A \emph{basis} of~\eqref{eq:lp} is a set~$B\subseteq[n]$ of~$m$ indices such that~$A_B$ is nonsingular.
The indices in~$B$ are called \emph{basic}, the other indices are \emph{nonbasic}.
Given a basis~$B$, the \emph{basic solution} is the unique $\ivec{x}{B}\in\Q^n$ with $\ivec{x}{B}_B=A_B^{-1}\vec{b}$ and $\icoord{x}{B}{j}=0$ for~$j\notin B$.
The basis is \emph{feasible} if $\ivec{x}{B}\geq\zero$, and a feasible basis is \emph{nondegenerate} if $\ivec{x}{B}_B>\zero$.
Its \emph{objective value} is $z_B\coloneqq\vec{c}^\top\ivec{x}{B}$, and a feasible basis is \emph{optimal} if $\ivec{x}{B}$ is an optimal solution of~\eqref{eq:lp}.

\begin{definition}\label{def:basis}
Given a basis~$B$ of~\eqref{eq:lp} and an index~$j\in[n]$,
\begin{enumerate}[label=(\alph*),ref=\thedefinition(\alph*)]
    \item\label{def:basis_a} the \emph{dual vector} is $\ivec{y}{B}\coloneqq A_B^{-\top}\vec{c_B}\in\Q^m$,
    \item\label{def:basis_b} the \emph{reduced cost} of index~$j$ is $\rc^B_j\coloneqq c_j-\vec{A_j}^\top\ivec{y}{B}$, 
    \item\label{def:basis_c} the index~$j$ is \emph{improving} if~$\rc^B_j<0$, and we write~$\Gamma(B)\coloneqq\{j\in[n]:\rc^B_j<0\}$,
    \item\label{def:basis_d} for nonbasic~$j$, the \emph{basic direction} is the vector $\ivec{d}{B,j}\in\Q^n$ with $\ivec{d}{B,j}_B=-A_B^{-1}\vec{A_j}$ and $\icoord{d}{B,j}{j}=1$, and with $\icoord{d}{B,j}{k}=0$ for all other nonbasic indices~$k\not\in B\cup\{j\}$,
    \item\label{def:basis_e} for feasible~$B$ and nonbasic improving~$j$, the \emph{ratio test} computes the step length
\begin{equation}
        \theta^{B,j}\coloneqq\max\{r\geq0:\ivec{x}{B}+r\ivec{d}{B,j}\geq\zero\}
        =\min\left\{\frac{\icoord{x}{B}{i}}{-\icoord{d}{B,j}{i}}\colon i\in B,\,\icoord{d}{B,j}{i}<0\right\}.
        \label{eq:ratio}
    \end{equation}
\end{enumerate}
\end{definition}

Throughout, we consider only \emph{nondegenerate} LPs, in which every feasible basis is nondegenerate, and we assume that the feasible region~$\polytope\coloneqq\{\vec{x}\in\R^n:\ A\vec{x}=\vec{b},\ \vec{x}\geq\zero\}$ is nonempty and bounded.
Then every linear function attains its minimum and its maximum over~$\polytope$ at basic solutions of feasible bases, and~$B$ is optimal if and only if~$\Gamma(B)=\emptyset$.
Further, given a feasible basis~$B$ and an improving index~$j$, exactly one index~$i^*\in B$ attains the minimum in the ratio test~\eqref{eq:ratio}, and~$0<\theta^{B,j}<\infty$.
The set~$B'\coloneqq (B\setminus\{i^*\})\cup\{j\}$ is a feasible basis, and
\begin{equation}
    \ivec{x}{B'}=\ivec{x}{B}+\theta^{B,j}\ivec{d}{B,j},\qquad
    z_{B'}=z_B+\theta^{B,j}\rc^B_j<z_B.
    \label{eq:pivot}
\end{equation}
We say that~$B'$ arises from~$B$ by the \emph{pivot} on~$j$, in which~$j$ \emph{enters} and~$i^*$ \emph{leaves} the basis.

Given a nondegenerate bounded linear program and a feasible basis, the simplex method iteratively performs a pivot on an improving index~$j\in\Gamma(B)$ of the current basis~$B$, and moves to the arising basis~$B'$ and its basic solution $\ivec{x}{B'}$, until reaching a basis without improving indices; see the pseudocode in Algorithm~\ref{alg:simplex}.%
\footnote{Algorithm~\ref{alg:simplex} is written so that its lines match those of Wolfe's method (Algorithm~\ref{alg:wolfe}).
The gray parts have no effect: the inner loop runs exactly once, and Line~\ref{line:sx-accept} does not change~$\vec{x}$.}
A \emph{simplex run} from a feasible basis~$B_0$ is a sequence~$B_0,B_1,\dots,B_f$ of feasible bases such that~$B_k$ arises from~$B_{k-1}$ by the pivot on some~$j_k\in\Gamma(B_{k-1})$ for every~$k\in[f]$, and~$\Gamma(B_f)=\emptyset$; its \emph{length} is the number~$f$ of pivots.
By~\eqref{eq:pivot}, every pivot strictly decreases the objective value, so no basis occurs twice, and the simplex method always terminates at an optimal basis.
The number of iterations depends on the \emph{pivot rule} that chooses the entering index in each step (Line~\ref{line:sx-choose}), in case there is more than one; see Table~\ref{tab:wolfe-rules} for some classical pivot rules and known lower bounds.
We refer to~\cite{BT97} for more details on the simplex method.

\begin{table}[tbp]
\centering
\caption{Insertion rules for Wolfe's method that are equivalent to classical pivot rules (Corollary~\ref{cor:intro-rules}), the atom~$j\in\Gamma(\vec{w_S})$ they insert at the corral~$S=S_k$ of a Wolfe run~$S_0,\dots,S_f$, and the resulting lower bounds on the number of iterations, in terms of the number~$n$ of variables of the LP.
Here~$j_\ell$ is the index inserted at the corral~$S_{\ell-1}$ for~$\ell\in[k]$, so~$j_1,\dots,j_k$ are the indices inserted before~$j$ ($j_0$ is a fixed initial index), and~$\pi$ is a uniformly random permutation of~$[n+1]$, drawn once.
}
\label{tab:wolfe-rules}
\small
\renewcommand{\arraystretch}{1.15}
\begin{tabularx}{\textwidth}{@{}l@{\hspace{2\tabcolsep}\quad}>{\raggedright\arraybackslash}Xll@{}}
\toprule
Insertion rule & Inserted atom~$j\in\Gamma(\vec{w_S})$ & Pivot rule & Lower bound\\
\midrule
linopt~\cite{Wolfe76} &~$j=\argmin\{\vec{w_S}^\top\ivec{p}{i}\colon i\in\Gamma(\vec{w_S})\}$ & Dantzig & $\Omega(2^{n/2})$~\cite{KM72}\\
least-index &~$j=\min\Gamma(\vec{w_S})$ & Bland & $\Omega(2^{n/2})$~\cite{AC78,AZ98}\\
least inserted &~$j=\argmin\{|\{\ell\in[k]\colon j_\ell=i\}|\colon i\in\Gamma(\vec{w_S})\}$ & Zadeh & $2^{\Omega(n)}$~\cite{DFH23}\\
round-robin &~$j=\argmin\{(i-j_k)\bmod(n+1)\colon i\in\Gamma(\vec{w_S})\}$ & Cunningham & $2^{\Omega(n)}$~\cite{AF17}\\
random &~$j\in\Gamma(\vec{w_S})$ uniformly at random & Random-Edge & $2^{\tilde\Omega(n^{1/4})}$~\cite{FHZ11}\\
random-order &~$j=\argmin\{\pi(i)\colon i\in\Gamma(\vec{w_S})\}$ & Random-Bland & $2^{\tilde\Omega(\sqrt{n})}$~\cite{FHZ14}\\
\bottomrule
\end{tabularx}
\end{table}

\subsection{Wolfe's minimum-norm-point method}
\label{sec:wolfe}

Wolfe's method~\cite{Wolfe74,Wolfe76} solves the minimum-norm-point problem~\eqref{eq:mnp} for given points~$\ivec{p}{1},\dots,\ivec{p}{n}\in\Q^D$, called \emph{atoms}.
For~$S\subseteq[n]$, we write~$\conv(S)$ and~$\aff(S)$ for the convex and the affine hull of the atoms~$\ivec{p}{i}$ with~$i\in S$.
A \emph{support} is a nonempty set~$S\subseteq[n]$, and it is \emph{affinely independent} if its atoms are.
The \emph{affine minimizer} of an affinely independent support~$S$ is the unique point~$\vec{w_S}\coloneqq\argmin_{\vec x\in\aff(S)}\norm{\vec x}$.
Its \emph{coefficients} $\ivec{\bar{\lambda}}{S}\in\R^S$ are the unique affine coefficients such that $\vec{w_S}=\sum_{i\in S}\icoord{\bar{\lambda}}{S}{i}\ivec{p}{i}$ and~$\sum_{i\in S}\icoord{\bar{\lambda}}{S}{i}=1$.
The affine minimizer and its coefficients can be computed in strongly polynomial time~\cite[Section~3.3]{Sch86}.
The support~$S$ is a \emph{corral} if $\icoord{\bar{\lambda}}{S}{i}>0$ for every~$i\in S$, that is, if~$\vec{w_S}$ lies in the relative interior of~$\conv(S)$.
Every singleton is a corral.

\begin{definition}\label{def:gap}
    Consider an~\eqref{eq:mnp} with atoms $\ivec{p}{1},\dots,\ivec{p}{n}\in\Q^D$.
The \emph{gap} of atom~$\ivec{p}{j}$ at a point~$\vec{w}\in\R^D$ is $g_j(\vec{w})\coloneqq\vec{w}^\top\ivec{p}{j}-\norm{\vec{w}}^2=\vec{w}^\top(\ivec{p}{j}-\vec{w})$.
The atom $\ivec{p}{j}$ is \emph{improving} at~$\vec{w}$ if~$g_j(\vec{w})<0$, and we write~$\Gamma(\vec{w})\coloneqq\{j\in[n]:g_j(\vec{w})<0\}$.
\end{definition}

Parts~(a) and~(c) of the following lemma are Wolfe's criteria for the affine and the convex hull~\cite[Lemmas~1.2 and~1.1]{DHR}.
\begin{lemma}\leavevmode\hskip-\labelsep\footnote{\label{fn:proofs}All missing proofs are deferred to the appendix.}\ %
\label{lem:affine-min}
    Consider~\eqref{eq:mnp} with atoms $\ivec{p}{1},\dots,\ivec{p}{n}\in\Q^D$.
Let~$S\subseteq[n]$ be an affinely independent support with affine minimizer~$\vec{w_S}$.
\begin{enumerate}
    \item[(a)] A point~$\vec{w}\in\aff(S)$ equals~$\vec{w_S}$ if and only if~$g_i(\vec{w})=0$ for every~$i\in S$.
    \item[(b)] If~$j\in[n]$ satisfies~$g_j(\vec{w_S})\neq0$, then~$S\cup\{j\}$ is affinely independent.
    \item[(c)] If~$S$ is a corral, then~$\vec{w_S}$ is the optimal solution of~\eqref{eq:mnp} if and only if~$\Gamma(\vec{w_S})=\emptyset$.
\end{enumerate}
\end{lemma}

Wolfe's method starts at a corral~$S_0$, for example at a single atom (cf.~Fig.~\ref{fig:wolfe-example} along with the following).
At a corral~$S$, the method stops if~$\Gamma(\vec{w_S})=\emptyset$, in which case~$\vec{w_S}$ is optimal by Lemma~\ref{lem:affine-min}(c).
Otherwise, a \emph{major cycle} adds an improving atom~$j$ to~$S$ (Line~\ref{line:wf-insert}), where~$j\notin S$ by Lemma~\ref{lem:affine-min}(a).
The resulting support~$S\cup\{j\}$ is affinely independent by Lemma~\ref{lem:affine-min}(b), but it is not necessarily a corral anymore.
As long as the support is not a corral, a \emph{minor cycle} (Lines~\ref{line:wf-direction}--\ref{line:wf-delete}) moves the coefficients towards those of the affine minimizer as far as nonnegativity permits and deletes one index whose coefficient has become zero; 
such an index exists, as some~$i\in S$ has $\icoord{\bar{\lambda}}{S}{i}\leq0$.
So every minor cycle deletes exactly one index.
Since singletons are corrals and since every affinely independent support contains at most~$D+1$ indices, the minor cycles end at a corral after at most~$D$ iterations.
At a corral, the current point moves to its affine minimizer (Line~\ref{line:wf-accept}), and the major cycle ends.
Consider Algorithm~\ref{alg:wolfe} for a pseudocode; it is the pseudocode of De~Loera et al.~\cite[Section~1.1]{DHR}, written in terms of the coefficient vector~$\vec{\lambda}$ of the current point and started at an arbitrary corral instead of a single atom.
We refer to~\cite[Section~1.1]{DHR} for more details on Wolfe's method.

A \emph{Wolfe run} from a corral~$S_0$ is a sequence~$S_0,S_1,\dots,S_f$ of corrals such that~$S_k$ is the corral at the end of a major cycle at~$S_{k-1}$ for every~$k\in[f]$, and~$\Gamma(\vec{w_{S_f}})=\emptyset$; its \emph{length} is the number~$f$ of major cycles.
Every major cycle strictly decreases the norm of the affine minimizer of the current corral~\cite{Wolfe76}, see~\cite[Section~1.1]{DHR}, so no corral occurs twice, and Algorithm~\ref{alg:wolfe} always terminates.
The number of major cycles depends on the \emph{insertion rule}, which chooses the improving atom in Line~\ref{line:wf-choose} (see Table~\ref{tab:wolfe-rules}).

\begin{figure}[t]
  \centering\scriptsize
  \renewcommand\mycommfont[1]{\scriptsize\ttfamily\textcolor{gray}{#1}}%
  \setlength{\algomargin}{2pt}%
  \begin{minipage}[t]{.497\linewidth}
    \begin{algorithm}[H]
    \begin{minipage}{\dimexpr\textwidth-13pt}
        \setstretch{1.0}

	    \caption{\algcapstrut\textsc{Simplex}($A$,~$\vec{b}$,~$\vec{c}$,~$B_0$)}\label{alg:simplex}
	    \Input{\algstrut{}nondegenerate~\eqref{eq:lp}, feasible basis~$B_0$}
	    \Output{\algstrut{}optimal basis~$B$, optimal solution~$\vec{x}=\ivec{x}{B}$}
	    \vgap
	    \hrule
	    \vgap
	    \setcounter{AlgoLine}{0}

	    \nl \algstrut$B\gets B_0$, $\vec{x}\gets\ivec{x}{B_0}$
        
        \SetInd{4pt}{4pt}
	    \nl \algstrut\While{$\Gamma(B)\neq\emptyset$}
	    {
	    \vspace{2pt}
		 \nl \algstrut{}choose~$j\in\Gamma(B)$\label{line:sx-choose}\tcp*{reduced cost~$\smash[b]{\rc^B_j}<0$}
		 \nl \algstrut$B\gets B\cup\{j\}$\label{line:sx-insert}\tcp*{insertion,~$x_j=0$}
		 \nl \algstrut\textcolor{grayif}{\While(\tcp*[f]{runs exactly once}){$B$ is not a basis\label{line:sx-loop}}
		 {
		 \nl \algstrut\textcolor{black}{$\vec{d}\gets\ivec{d}{B\setminus\{j\},j}$}\label{line:sx-direction}\tcp*{basic direction}
		 \nl \algstrut\textcolor{black}{$\theta\gets\max\{r\geq0:\vec{x}+r\vec{d}\geq\zero\}$}\label{line:sx-step}\tcp*{ratio test}
		 \nl \algstrut\textcolor{black}{$\vec{x}\gets\vec{x}+\theta\vec{d}$}\label{line:sx-move}\;
		 \nl \algstrut\textcolor{black}{$B\gets B\setminus\{i\}$ for an~$i\in B$ with~$x_i=0$}\label{line:sx-delete}\;
		 }}
		 \nl \algstrut\textcolor{grayif}{$\vec{x}\gets\ivec{x}{B}$}\label{line:sx-accept}\tcp*{no change}
	    }
	    \nl \algstrut\Return~$B$,~$\vec{x}$
	    \vspace{2pt}
    \end{minipage}
	\end{algorithm}
  \end{minipage}%
  \hfill
  \begin{minipage}[t]{.497\linewidth}
    \begin{algorithm}[H]
    \begin{minipage}{\dimexpr\textwidth-13pt}
        \setstretch{1.0}

	    \caption{\algcapstrut\textsc{Wolfe}($\ivec{p}{1}\!,\dots\!,\ivec{p}{n}\!,\!S_0$)}\label{alg:wolfe}
	    \Input{\algstrut{}atoms~$\ivec{p}{1},\dots,\ivec{p}{n}\in\Q^D$, corral~$S_0$}
	    \Output{\algstrut{}corral~$S$, optimal solution~$\vec{w_S}$}
	    \vgap
	    \hrule
	    \vgap
	    \setcounter{AlgoLine}{0}

	    \nl \algstrut$S\gets S_0$, $\vec{\lambda}\gets\ivec{\bar{\lambda}}{S_0}$
        
        \SetInd{4pt}{4pt}
	    \nl \algstrut\While(\tcp*[f]{major cycles}){$\Gamma(\vec{w_S})\neq\emptyset$}
	    {
	    \vspace{2pt}
		 \nl \algstrut{}choose~$j\in\Gamma(\vec{w_S})$\label{line:wf-choose}\tcp*{gap~$g_j(\vec{w_S})<0$}
		 \nl \algstrut$S\gets S\cup\{j\}$\label{line:wf-insert}\tcp*{insertion,~$\lambda_j=0$}
		 \nl \algstrut\While(\tcp*[f]{minor cycles}){$S$ is not a corral\label{line:wf-loop}}
		 {
		 \nl \algstrut$\vec{d}\gets\ivec{\bar{\lambda}}{S}-\vec{\lambda}$\label{line:wf-direction}\tcp*{to affine minimizer}
		 \nl \algstrut$t\gets\max\{r\in[0,1]:\vec{\lambda}+r\vec{d}\geq\zero\}$\label{line:wf-step}
         
		 \nl \algstrut$\vec{\lambda}\gets\vec{\lambda}+t\vec{d}$\label{line:wf-move}\;
		 \nl \algstrut$S\gets S\setminus\{i\}$ for an~$i\in S$ with~$\lambda_i=0$\label{line:wf-delete}\;
		 }
		 \nl \algstrut$\vec{\lambda}\gets\ivec{\bar{\lambda}}{S}$\label{line:wf-accept}\tcp*{acceptance}
	    }
	    \nl \algstrut\Return~$S$,~$\vec{w_S}$
	    \vspace{2pt}
    \end{minipage}
	\end{algorithm}
  \end{minipage}
\end{figure}

\section{Wolfe's method generalizes the simplex method}
\label{sec:equivalence}

This section shows how to transform a linear program of the form~\eqref{eq:lp} into an instance of the minimum-norm-point problem of the form~\eqref{eq:mnp} such that there is a one-to-one correspondence between simplex executions for the former and executions of Wolfe's method for the latter.
We require the input linear program to be \emph{admissible}, that is
\begin{enumerate}[label=(A\arabic*),ref=(A\arabic*),leftmargin=*]
    \item\label{ass:integral}
    its data are integral, i.e.,~$A\in\Z^{m\times n}$,~$\vec{b}\in\Z^m$, and~$\vec{c}\in\Z^n$;
    \item\label{ass:bounded}
    its feasible region~$\polytope=\{\vec{x}\in\R^n:A\vec{x}=\vec{b},\
    \vec{x}\geq\zero\}$ is nonempty and bounded;
    \item\label{ass:nondeg}
    all its feasible bases are nondegenerate.
\end{enumerate}
For rational input, the requirement~\ref{ass:integral} can be achieved by scaling the constraint rows and the objective appropriately.
The requirement~\ref{ass:nondeg} is also unproblematic, in the sense that all the worst-case constructions of linear programs that we aim to carry over to MNPs satisfy it. 
Some of these constructions do not satisfy \ref{ass:bounded}, as their feasible sets are unbounded, which we resolve by adding a new variable with a constraint that bounds the feasible region without changing the simplex runs (Appendix \ref{app:cor-rules}).

Our construction proceeds in two parts.
First, using standard LP transformations, we reformulate~\eqref{eq:lp} as (Section~\ref{sec:hom})
\begin{equation}
    \textstyle
    \min\ \vec{\gamma}^\top\vec{\lambda},\qquad\text{subject to}\quad
    \sum_j\lambda_j\,\ivec{\alpha}{j}=\zero,\quad
    \sum_j\lambda_j=1,\quad\vec{\lambda}\geq\zero,
    \label{eq:induced-lp}
    \tag{LP$_\lambda$}
\end{equation}
with~$\vec{\gamma} > \vec{0}$ and~$\ivec{\alpha}{j}\in\Q^m$.
Note that \ref{eq:induced-lp} more closely resembles~\eqref{eq:mnp} (and~\eqref{eq:CQM}): Its constraints demand that the variables~$\vec{\lambda}$ yield a convex combination of the origin using the vectors~$\ivec{\alpha}{j}$.
We obtain this form in three steps:
\begin{enumerate}
    \item Ensure~$\sum_j \lambda_j = 1$ and~$\vec{\lambda} \geq \vec{0}$:
    Choose~$T>\one^\top\vec{x}$ for every feasible~$\vec{x}$, set~$\lambda_j \coloneqq x_j/T$ for~$j \in \{1,\dots,n\}$ and introduce
   ~$\lambda_{n+1} \coloneqq 1-\sum_{j=1}^n\lambda_j$.
    
    \item Ensure right-hand side is~$\vec{0}$:
    Using~$\sum_j\lambda_j=1$, rewrite
   ~$A\vec{x}=\vec{b}$ as
   ~$\sum_j\lambda_j\,\ivec{\alpha}{j}=\zero$, where
   ~$\ivec{\alpha}{j} \coloneqq \vec{A_j}-\vec{b}/T$ for~$j\le n$
    and~$\ivec{\alpha}{n+1} \coloneqq -\vec{b}/T$.
    
    \item Ensure~$\vec{\gamma} > \vec{0}$:
    Set~$\gamma_j \coloneqq c_j + C$ for~$j \in \{1,\dots,n\}$ and
   ~$\gamma_{n+1} \coloneqq C$ for some~$C > \max_{j\in[n]}|c_j|$.
    This adds the constant~$\sum_j C\lambda_j = C$ to the objective.
\end{enumerate}
The key part of the analysis is the proof that these transformations preserve simplex executions (Lemma~\ref{lem:hom}).

The second part of our construction is the encoding of~\eqref{eq:induced-lp} as an instance of the minimum-norm-point problem of the form~\eqref{eq:mnp} (Section~\ref{sec:atoms}).
For each variable~$\lambda_j$, we introduce an atom
$\ivec{p}{j} \coloneqq \smash{\left(\begin{smallmatrix}N\gamma_j\\N^2\ivec{\alpha}{j}\end{smallmatrix}\right)} \in \mathds{Q}^{m+1}$, such that the Pythagorean theorem yields
\[
\textstyle
\bigl\lVert\sum_j\lambda_j\,\ivec{p}{j}\bigr\rVert^2
=
N^2(\vec{\gamma}^\top\vec{\lambda})^2
+
N^4\bigl\lVert\sum_j
\lambda_j\,\ivec{\alpha}{j}\bigr\rVert^2.
\]
Since~$\vec{\gamma} > 0$, minimizing the first term is equivalent to minimizing the objective of~\eqref{eq:induced-lp}.
The second term measures the violation of the equality constraints of~\eqref{eq:induced-lp}. 
By increasing~$N$ we can make these violations more expensive in relation to the objective.
For the atoms corresponding to any feasible basis, the coefficients
of their affine minimizer then approach the LP's basic feasible
solution. 
We choose~$N$ large enough so that each simplex pivot can be
reproduced by Wolfe's method: the atom corresponding to the entering
variable is inserted, and the atom corresponding to the leaving
variable is deleted (see Fig.~\ref{fig:simplex-wolfe}).
We determine a sufficiently large~$N$ of polynomial encoding length
via standard determinant bounds using Cramer's rule and Hadamard's
inequality (Lemmas~\ref{lem:bases-corrals} and~\ref{lem:encoding}).

\begin{figure}[tbp]
    \centering
    \begingroup
    \definecolor{edgecol}{gray}{0.55}
    \resizebox{\linewidth}{!}{%
\begin{tikzpicture}[
    every node/.style={font=\small},
    axis/.style={black!20,line width=.5pt,-{Stealth[length=1.5mm]}},
    axislabel/.style={text=black!45,inner sep=1pt},
    trajectory/.style={wolfecol,line width=1.05pt,
        -{Stealth[length=1.6mm,width=1.2mm]},
        shorten <=2pt,shorten >=2.2pt},
    wolfetrajectory/.style={wolfecol,line width=1.05pt,
        -{Stealth[length=.9mm,width=.75mm]},
        shorten <=1.8pt,shorten >=1.8pt}]
  \path[use as bounding box] (0,-.2) rectangle (13.8,4.7);

  \begin{scope}[shift={(.95,1.25)},x=2cm,y=2cm]
    \coordinate (v3) at (0,0);
    \coordinate (v1) at (2,0);
    \coordinate (v2) at (0,1);
    \draw[axis] (-.1,0)--(2.35,0) node[right,axislabel] {$x_1$};
    \draw[axis] (0,-.1)--(0,1.25) node[above,axislabel] {$x_2$};
    \fill[black!3] (v3)--(v1)--(v2)--cycle;
    \draw[edgecol,line width=.7pt] (v3)--(v1)--(v2)--cycle;
    \draw[trajectory] (v3)--(v2);
    \draw[trajectory] (v2)--(v1);
    \foreach \v in {v3,v2}\fill[wolfecol] (\v) circle[radius=2pt];
    \draw[wolfecol,line width=.8pt,fill=white] (v1) circle[radius=2.5pt];
    \fill[wolfecol] (v1) circle[radius=1.3pt];
    \node[below left=2pt,inner sep=1pt] at (v3) {$\left(\begin{smallmatrix}0\\0\end{smallmatrix}\right)$};
    \node[below right=2pt,inner sep=1pt] at (v1) {$\left(\begin{smallmatrix}2\\0\end{smallmatrix}\right)$};
    \node[left=4pt,inner sep=1pt] at (v2) {$\left(\begin{smallmatrix}0\\1\end{smallmatrix}\right)$};
  \end{scope}

  \begin{scope}[shift={(8.8,1.28)},x=.21cm,y=.21cm]
    \coordinate (O) at (0,0);
    \coordinate (p1) at (6,{27/5});
    \coordinate (p2) at (3,{72/5});
    \coordinate (p3) at (12,{27/5});
    \coordinate (p4) at (12,{-18/5});
    \coordinate (w0) at (12,0);
    \coordinate (z1) at ({51/5},0);
    \coordinate (w1) at ({204/25},{102/25});
    \coordinate (z2) at ({36/5},{18/5});
    \coordinate (w2) at ({432/65},{288/65});
    \draw[axis] (-.6,0)--(13.2,0) node[right,axislabel] {$x_1$};
    \draw[axis] (0,-4.2)--(0,14.5) node[above,axislabel] {$x_2$};
    \fill[black!3] (p2)--(p3)--(p4)--(p1)--cycle;
    \draw[edgecol,line width=.7pt] (p2)--(p3)--(p4)--(p1)--cycle;
    \draw[edgecol,line width=.6pt,dash pattern=on 2pt off 2.5pt] (p2)--(p4);
    \foreach \a/\b in {w0/z1,z1/w1,w1/z2,z2/w2}
      \draw[wolfetrajectory] (\a)--(\b);
    \foreach \p in {p1,p2,p3,p4}\fill (\p) circle[radius=2pt];
    \foreach \p in {w0,w1}\fill[wolfecol] (\p) circle[radius=1.7pt];
    \foreach \p in {z1,z2}\draw[wolfecol,line width=.7pt,fill=white] (\p) circle[radius=1.4pt];
    \draw[wolfecol,line width=.7pt,fill=white] (w2) circle[radius=2.1pt];
    \fill[wolfecol] (w2) circle[radius=1pt];
    \fill (O) circle[radius=1.5pt];
    \node[above left=2pt,xshift=5pt,inner sep=1pt] at (p2) {$\ivec{p}{2}$};
    \node[left=3pt,inner sep=1pt] at (p1) {$\ivec{p}{1}$};
    \node[right=3pt,inner sep=1pt] at (p3) {$\ivec{p}{3}$};
    \node[right=3pt,inner sep=1pt] at (p4) {$\ivec{p}{4}$};
    \node[above left=1pt,inner sep=1pt] at (O) {$\vec{0}$};
  \end{scope}

  \node at (3.15,.02) {bases:~$\{3\}\to\{2\}\to\{1\}$};
  \node at (10.55,.02) {corrals:~$\{3,4\}\to\{2,4\}\to\{1,4\}$};
\end{tikzpicture}
    }
    \endgroup
    \setlength{\abovecaptionskip}{2pt}
    \caption{Execution of the simplex method with Dantzig's rule (left) for the linear program~$\min\{-2x_1-3x_2:x_1+2x_2\le2,\ x_1,x_2\ge0\}$; and execution of Wolfe's method with the linopt rule (right) for the corresponding MNP with atoms~$\ivec{p}{1}=\smash{\tfrac15\left(\begin{smallmatrix}30\\27\end{smallmatrix}\right)}$,
$\ivec{p}{2}=\smash{\tfrac15\left(\begin{smallmatrix}15\\72\end{smallmatrix}\right)}$,~$\ivec{p}{3}=\smash{\tfrac15\left(\begin{smallmatrix}60\\27\end{smallmatrix}\right)}$,
and~$\ivec{p}{4}=\smash{\tfrac15\left(\begin{smallmatrix}60\\-18\end{smallmatrix}\right)}$.
For illustration we use~$T=5$,~$C=4$, and scale~$N=3$ (faithful scale would be~$N\approx44 \cdot 10^{12}$).
In Wolfe's method, filled blue circles mark the positions at which indices
$2$ and~$1$ are inserted into the support, in this order.
At the open circles, indices~$3$ and~$2$ are deleted.}
    \label{fig:simplex-wolfe}
\end{figure}

The key difficulty is to establish the exact correspondence of
executions despite the small constraint violations:
feasible bases yield corrals, and each simplex pivot corresponds
to one Wolfe major cycle with the same entering and leaving
indices, and conversely, from corresponding initial bases and
corrals (Theorem~\ref{thm:equivalence}).
Moreover, the gaps are the reduced costs multiplied by a common
positive factor.
Together with the preserved indices and history, this translates
pivot rules invariant under such scaling into insertion rules;
Lemma~\ref{lem:rules} makes this correspondence explicit for the
rules of Table~\ref{tab:wolfe-rules}.

\subsection{The induced linear program}
\label{sec:hom}

We now bring the linear program into the form that resembles \eqref{eq:mnp} more closely (see the discussion above).
In the rest of this section, \eqref{eq:lp} denotes a fixed admissible linear program,
and we use the notation from \ref{ass:integral}--\ref{ass:nondeg}.
Let~$U\geq1$ be the largest absolute value of an entry of~$A$ or~$\vec{b}$.
It is well-known that every feasible point~$\vec{x}\in\polytope$ satisfies
\begin{equation}
    2\cdot\one^\top\vec{x}<2m(mU)^m+1\eqqcolon T,
    \label{eq:T}
\end{equation}
but to be self-contained we give the proof in Lemma~\ref{lem:T-bound} in Appendix~\ref{app:statements}.

Similar to \cite[Lemmas~2.6 and~2.8]{DHR}, we construct a linear program that is equivalent to the given admissible program by four textbook operations:
adding a new variable~$x_{n+1}\geq0$ and the constraint $\one^\top\vec{x}+x_{n+1}=T$, whose slack~$x_{n+1}$ is positive on~$\polytope$ by~\eqref{eq:T}; substituting $(\vec{x}^\top,x_{n+1})^\top=T\vec{\lambda}$ and dividing the constraints and the objective by~$T$; subtracting~$\vec{b}/T$ times the new constraint from the others; and adding~$C$ times it to the objective, for the constant~$C$ in~\eqref{eq:gamma} below.
The last step makes all objective coefficients positive, so that minimizing the objective value is the same as minimizing its square, which we will use in Section~\ref{sec:atoms}.
This results in the \emph{induced linear program}~\eqref{eq:induced-lp},
where the vectors $\ivec{\alpha}{j}\in\Q^m$ and the objective coefficients~$\gamma_j\in\Z_{\geq1}$ are given by
\begin{align}
    \ivec{\alpha}{n+1}&=-\frac{\vec{b}}{T},&
    \ivec{\alpha}{j}&=\vec{A_j}-\frac{\vec{b}}{T}\quad(j\in[n]),
    \label{eq:alpha}\\
    \gamma_{n+1}&=1+\max_{j\in[n]}|c_j|\eqqcolon C,&
    \gamma_j&=C+c_j\quad(j\in[n]).
    \label{eq:gamma}
\end{align}
Of course, we can also write~\eqref{eq:induced-lp} in standard form $\min\{\vec{\gamma}^\top\vec{\lambda}:M\vec{\lambda}=(\zero^\top,1)^\top,\ \vec{\lambda}\geq\zero\}$, where~$M\in\Q^{(m+1)\times(n+1)}$ has columns $\vec{M_j}=((\ivec{\alpha}{j})^\top,1)^\top$,~$j\in[n+1]$.
Then, all notions of Section~\ref{sec:simplex} apply to every basis~$S\subseteq[n+1]$ of~\eqref{eq:induced-lp}, with~$A$,~$\vec{b}$,~$\vec{c}$,~$[n]$ replaced by the data~$M$, $(\zero^\top,1)^\top$,~$\vec{\gamma}$,~$[n+1]$.
The resulting quantities can be found in Table~\ref{tab:basis_identities}, which introduces the notation that we will use from now on.

\begin{table}[tbp]
\centering
\caption{The notions of Section~\ref{sec:simplex} applied to~\eqref{eq:lp} and~\eqref{eq:induced-lp}.
}
\label{tab:basis_identities}
\small
\renewcommand{\arraystretch}{1.2}
\begin{tabular*}{\textwidth}{@{\extracolsep{\fill}}l>{$\displaystyle}l<{$}>{$\displaystyle}l<{$}r@{}}
    \toprule
    & \text{basis~$B$ of~\eqref{eq:lp}} & \text{basis~$S$ of~\eqref{eq:induced-lp}} & \\
    \midrule
    basic solution & \ivec{x}{B}_B=A_B^{-1}\vec{b} & \ivec{\lambda}{S}_S=M_S^{-1}(\zero^\top,1)^\top & \tabeq{eq:nf-solution}\\
    dual vector & \ivec{y}{B}=A_B^{-\top}\vec{c_B} & ((\ivec{v}{S})^\top,h^S)^\top=M_S^{-\top}\vec{\gamma_S} & \tabeq{eq:nf-dual}\\
    reduced cost & \rc^B_j=c_j-\vec{A_j}^\top\ivec{y}{B} & \rg^S_j=\gamma_j-(\ivec{v}{S})^\top\ivec{\alpha}{j}-h^S & \tabeq{eq:nf-rc}\\
    basic direction & \ivec{d}{B,j}_B=-A_B^{-1}\vec{A_j} & \ivec{d}{S,j}_S=-M_S^{-1}\vec{M_j} & \tabeq{eq:nf-direction}\\
    ratio test & \theta^{B,j}=\min_{i\in B:\,\icoord{d}{B,j}{i}<0}\left\{{\icoord{x}{B}{i}}/{-\icoord{d}{B,j}{i}}\right\} & \theta^{S,j}=\min_{i\in S:\,\icoord{d}{S,j}{i}<0}\left\{{\icoord{\lambda}{S}{i}}/{-\icoord{d}{S,j}{i}}\right\} & \tabeq{eq:ratio-induced}\\[2pt]
    \bottomrule
\end{tabular*}
\end{table}

For~$B\subseteq[n]$, we write~$\hat{B}\coloneqq B\cup\{n+1\}$, and call~$n+1$ the \emph{slack index}.
The simplex runs on~\eqref{eq:lp} and on~\eqref{eq:induced-lp} coincide:~$B$ is a feasible basis of~\eqref{eq:lp} if and only if~$\hat{B}$ is one of~\eqref{eq:induced-lp}, nonbasic indices have the same reduced costs, and the leaving indices coincide for every pivot.
This is formalized in the following lemma.

\begin{lemma}\proofinappendix\label{lem:hom}
The induced program~\eqref{eq:induced-lp} is nondegenerate, its constraint matrix has full row rank, and its feasible set is nonempty and bounded.
The set~$B$ is a feasible basis of~\eqref{eq:lp} if and only if~$\hat{B}= B\cup\{n+1\}$ is a feasible basis of~\eqref{eq:induced-lp}, in which case
\begin{enumerate}
    \item[(a)] their basic solutions satisfy\quad $\ivec{\lambda}{\hat{B}}=\bigl((\ivec{x}{B})^\top/T,\,1-\one^\top\ivec{x}{B}/T\bigr)^\top$,
    \item[(b)] the dual vectors satisfy\quad $((\ivec{v}{\hat{B}})^\top,h^{\hat{B}})^\top=((\ivec{y}{B})^\top,\,C+z_B/T)^\top$,
    \item[(c)] the reduced costs satisfy\quad~$\rg^{\hat{B}}_j=\rc^B_j$ for~$j\in[n]$, and~$\rg^{\hat{B}}_{n+1}=0$.
\end{enumerate}
For every~$j\in\Gamma(B)$, the pivot on~$j$ at~$\hat{B}$ has the same unique leaving index as the one at~$B$.
Thus,~$B_0,\dots,B_f$ is a simplex run on~\eqref{eq:lp} if and only if~$\hat{B}_0,\dots,\hat{B}_f$ is one on~\eqref{eq:induced-lp}.
\end{lemma}

\subsection{The induced minimum-norm-point problem}
\label{sec:atoms}

Now we give the minimum-norm-point problem that is a penalty version of the linear program from the previous section (see the discussion at the start of this section).
The objective coefficients in the induced linear program~\eqref{eq:induced-lp} are integers~$\gamma_j\geq1$, and the vectors $T\ivec{\alpha}{j}$ are integral, for all~$j\in[n+1]$, by~\ref{ass:integral}, \eqref{eq:alpha}, and~\eqref{eq:gamma}.
For a scale factor~$N>0$, which is an integer multiple of~$T$ that we fix in Step~3 below, we define the atoms
\begin{equation}    \ivec{p}{j}=\bigl(N\gamma_j,\,N^2(\ivec{\alpha}{j})^\top\bigr)^\top\in\Z^{m+1},\qquad j\in[n+1],
    \label{eq:atoms}
\end{equation}
and give a name to the minimum-norm-point problem on these atoms
\begin{equation}
    \min\ \norm{\vec{w}},\quad\text{subject to}\quad\vec{w}\in\conv\{\ivec{p}{1},\dots,\ivec{p}{n+1}\}.
    \label{eq:mnp-lambda}
    \tag{MNP$_\lambda$}
\end{equation}

All notions of Section~\ref{sec:wolfe} apply to these atoms, with~$[n]$ replaced by~$[n+1]$.
For a point~$\vec{w}\in\R^{m+1}$, we say that~$w_1$ is its \emph{objective coordinate}, and~$w_2,\dots,w_{m+1}$ are its \emph{constraint coordinates} (one for each of the~$m$ rows of the constraint $\sum_j\lambda_j\ivec{\alpha}{j}=\zero$ of~\eqref{eq:induced-lp}).
For~$\vec{\lambda}\in\R^{n+1}$ with $\one^\top\vec{\lambda}=1$, we obtain the \emph{penalty identity}
\begin{equation}
    \begin{gathered}
        \sum_{j\in[n+1]}\lambda_j\ivec{p}{j}=\Bigl(N\vec{\gamma}^\top\vec{\lambda},\ N^2\sum_{j\in[n+1]}\lambda_j(\ivec{\alpha}{j})^\top\Bigr)^\top,\\
        \Bigl\lVert\sum_{j\in[n+1]}\lambda_j\ivec{p}{j}\Bigr\rVert^2
           =N^2(\vec{\gamma}^\top\vec{\lambda})^2
              +N^4\Bigl\lVert\sum_{j\in[n+1]}\lambda_j\ivec{\alpha}{j}\Bigr\rVert^2.
    \end{gathered}
    \label{eq:penalty}
\end{equation}
Thus, the minimum-norm-point problem \eqref{eq:mnp-lambda} is a penalty version of~\eqref{eq:induced-lp} with a scale factor~$N$.
We show, in four steps, that one major cycle of Wolfe's method on~\eqref{eq:mnp-lambda}
corresponds to one pivot of the simplex method on \eqref{eq:induced-lp}, already for a scale factor~$N$ of polynomial bit length.
Steps~1 and~2 hold for every~$N>0$.

\paragraph{Step 1: Reduced costs and gaps.}
Wolfe's method chooses the entering index by its gap, whereas the simplex method chooses it by its reduced cost (Line~\ref{line:sx-choose} of Algorithms~\ref{alg:simplex} \&~\ref{alg:wolfe}).
In this step, we prove a close connection between these two quantities.

Let~$S$ be a basis of~\eqref{eq:induced-lp}, and let~$j$ be a nonbasic index.
Since the basic direction~$\ivec{d}{S,j}$ satisfies the homogeneous constraints and has objective value $\vec{\gamma}^\top\ivec{d}{S,j}={}\rg^S_j$
(Lemma~\ref{lem:induced-facts}(c) in Appendix~\ref{app:statements}), the atom $\ivec{p}{j}$ is an affine combination of the atoms of~$S$, shifted by~$N\rg^S_j$ in the objective coordinate.

\begin{lemma}\proofinappendix\label{lem:direction}
Let~$S$ be a basis of~\eqref{eq:induced-lp} and~$j\in [n+1]\setminus S$.
Then, for every~$\vec{w}\in\R^{m+1}$,
\begin{equation}
    \ivec{p}{j}+\sum_{i\in S}\icoord{d}{S,j}{i}\ivec{p}{i}=N\rg^S_j\vec{e_1},\qquad
    g_j(\vec{w})+\sum_{i\in S}\icoord{d}{S,j}{i}g_i(\vec{w})=N\rg^S_jw_1.
    \label{eq:identity}
\end{equation}
\end{lemma}

Thus, at every point~$\vec{w}$ at which all atoms of~$S$ have gap zero, every index~$k\in[n+1]$ has the gap~$g_k(\vec{w})=N\rg^S_kw_1$ (for~$k\in S$, both sides vanish).
By Lemma~\ref{lem:affine-min}(a), the affine minimizer~$\vec{w_S}$ is such a point.\footnote{Lemma~\ref{lem:affine-min}(a) applies, since the atoms of~$S$ are affinely independent by Lemma~\ref{lem:corral}(a) below.}
In Step~2, we compute~$\vec{w_S}$.

\paragraph{Step 2: The affine minimizer of a basis.}
Let~$S$ be a basis of~\eqref{eq:induced-lp}.
Its atoms are affinely independent (Lemma~\ref{lem:corral}(a)) and, by the penalty identity~\eqref{eq:penalty}, the coefficients~$\vec{\lambda}\in\R^S$ of the affine minimizer~$\vec{w_S}$ minimize~$N^2(\vec{\gamma_S}^\top\vec{\lambda})^2+N^4\norm{\vec{\delta}}^2$, where $\vec{\delta}\coloneqq\sum_{i\in S}\lambda_i\ivec{\alpha}{i}$ is the \emph{violation} of the homogeneous constraints.
The basic solution $\ivec{\lambda}{S}$ has no violation and the objective value $\vec{\gamma}^\top\ivec{\lambda}{S}={}h^S$
(Lemma~\ref{lem:induced-facts}(a),(b) in Appendix~\ref{app:statements}), while the affine minimizer might use a small violation.
A violation~$\vec{\delta}$ changes the objective value by $\vec{\delta}^\top\ivec{v}{S}$ (by~\eqref{eq:dual} in Appendix~\ref{app:statements}), so the best violation is of the form $-\varepsilon\ivec{v}{S}$ for some~$\varepsilon\in\R$, which changes the coefficients by $-\varepsilon\ivec{u}{S}$ for the vector $\ivec{u}{S}$ in~\eqref{eq:projection-quantities} below.
Minimizing over~$\varepsilon$ gives the following optimal value~$\varepsilon^S$ and the objective value~$s^S$ of the affine minimizer
\begin{equation}
    \ivec{u}{S}=M_S^{-1}((\ivec{v}{S})^\top,0)^\top\in\Q^S,\qquad
    \varepsilon^S=\frac{h^S}{N^2+\norm{\ivec{v}{S}}^2},\qquad
    s^S=N^2\varepsilon^S.
    \label{eq:projection-quantities}
\end{equation}
We verify these quantities in the following lemma, which also shows that the objective coordinate of~$\vec{w_S}$ is positive if~$S$ is feasible.
Thus, by Step~1, the gaps at~$\vec{w_S}$ are the reduced costs times a common positive factor.

\begin{lemma}\proofinappendix\label{lem:corral}
Let~$S$ be a basis of~\eqref{eq:induced-lp}, and let $\ivec{u}{S}$,~$\varepsilon^S$, and~$s^S$ be given by~\eqref{eq:projection-quantities}.
Then,
\begin{enumerate}
    \item[(a)] the atoms $\ivec{p}{i}$,~$i\in S$, are affinely independent, and the affine minimizer~$\vec{w_S}$ and its coefficients $\ivec{\bar{\lambda}}{S}$ satisfy
\begin{equation}
        \begin{gathered}
        \ivec{\bar{\lambda}}{S}=\ivec{\lambda}{S}-\varepsilon^S\ivec{u}{S},~~
        \vec{w_S}=\left(\begin{smallmatrix}Ns^S\\-N^2\varepsilon^S\ivec{v}{S}\end{smallmatrix}\right),~~
        \vec{\gamma}^\top\ivec{\bar{\lambda}}{S}=s^S=h^S-\varepsilon^S\norm{\ivec{v}{S}}^2,
        \end{gathered}
        \label{eq:projection}
    \end{equation}
    \item[(b)] every index~$j\in[n+1]$ satisfies~$g_j(\vec{w_S})=N^2s^S\rg^S_j$; if~$S$ is feasible, then~$s^S>0$.
\end{enumerate}
\end{lemma}
In particular, at the affine minimizer~$\vec{w_S}$ of a feasible basis~$S$, the improving atoms are those whose indices are improving at the basis~$S$, and comparing gaps (as the linopt rule does) is comparing reduced costs (as Dantzig's rule does).

\paragraph{Step 3: Feasible bases are corrals.}
By Step~2, the coefficients of the affine minimizer of a basis~$S$ are the basic solution shifted by $(-\varepsilon^S\ivec{u}{S})$.
We now fix~$N$ so large that, for a feasible basis~$S$, this shift is far below the smallest basic coordinate.
The quantities in \eqref{eq:nf-solution}--\eqref{eq:ratio-induced} are computed from solutions of linear systems with the matrix~$M_S$ or its transpose.
After multiplying the first~$m$ rows of~$M$ by~$T$, these systems become integral, and their entries are entries of the integer vectors $(\gamma_j,T(\ivec{\alpha}{j})^\top,1)^\top$ with~$j\in[n+1]$.
Let~$K\geq1$ be the largest absolute value of these entries, and set
\begin{equation}
    H\coloneqq\bigl((m+1)K^2\bigr)^{m+1},\qquad N\coloneqq2mTH^3.
    \label{eq:scale}
\end{equation}
By Hadamard's inequality and Cramer's rule, every nonzero coordinate of a basic solution or of a basic direction of~\eqref{eq:induced-lp} has absolute value between~$1/H$ and~$H$ (proof of Lemma \ref{lem:bases-corrals}), so the basic coordinates of a feasible basis are at least~$1/H$ by nondegeneracy.
By the choice of~$N$, the shift $\varepsilon^S\ivec{u}{S}$ is far below~$1/H$ in every coordinate (Lemma \ref{lem:bases-corrals}(c)).
We choose~$N$ larger than needed here, because we need the smaller shift in Step~4.

\begin{lemma}\proofinappendix\label{lem:bases-corrals}
Let~$S$ be a feasible basis of~\eqref{eq:induced-lp}.
Then~$\icoord{\bar{\lambda}}{S}{i}>0$ for every~$i\in S$, so~$S$ is a corral.
More precisely, for every~$i\in S$ and~$j\in [n+1]\setminus S$,
\begin{enumerate}
    \item[(a)] we have $\icoord{\lambda}{S}{i}\geq1/H$,
    \item[(b)] the nonzero entries of $\ivec{d}{S,j}$ have absolute values between~$1/H$ and~$H$,
    \item[(c)] we have $\norm{\ivec{\bar{\lambda}}{S}-\ivec{\lambda}{S}}_\infty=\varepsilon^S\norm{\ivec{u}{S}}_\infty\leq1/(4mH^3)$.
\end{enumerate} 
\end{lemma}

\paragraph{Step 4: One major cycle is one pivot.}
Let~$S$ be a feasible basis of~\eqref{eq:induced-lp}, let~$\rg^S_j<0$, and let~$S'$ be the pivot on~$j$, with leaving index~$i^*$.
By Steps~2 and~3, $\ivec{p}{j}$ is improving at~$\vec{w_S}$, and~$S$ and~$S'$ are corrals.
After inserting $\ivec{p}{j}$, the support~$S\cup\{j\}$ consists of~$m+2$ affinely independent atoms in~$\R^{m+1}$, so its affine hull is~$\R^{m+1}$ and its affine minimizer is the origin.
The coefficients of the origin lie on the ray along which the simplex method moves, beyond the next basic solution $\ivec{\lambda}{S'}$ (proof of Lemma~\ref{lem:pivot}).
The minor cycle moves on a segment close to this ray, since the shift of Step~3 is tiny, so it also deletes~$i^*$.

\begin{lemma}\proofinappendix\label{lem:pivot}
Let~$S$ be a feasible basis of~\eqref{eq:induced-lp}, let~$\rg^S_j<0$, and let~$S'=(S\setminus\{i^*\})\cup\{j\}$ be the pivot on~$j$.
Then, the Wolfe major cycle that inserts $\ivec{p}{j}$ at the corral~$S$ consists of exactly one minor cycle which deletes~$i^*$.
The next corral is~$S'$.
\end{lemma}

The atoms in~\eqref{eq:atoms} are integer vectors of polynomial bit length, and they can be computed in polynomial time (Lemma~\ref{lem:encoding} in Appendix~\ref{app:statements}).
Lemmas~\ref{lem:hom}, \ref{lem:corral}, \ref{lem:bases-corrals}, and~\ref{lem:pivot} for the bases~$\hat{B}$ of~\eqref{eq:induced-lp} give the following detailed version of Theorem~\ref{thm:intro-equivalence}.

\setcounter{theorem}{0}%
\begin{theorem}\proofinappendix\label{thm:equivalence}
The atoms~$\ivec{p}{1},\ldots,\ivec{p}{n+1}\in\Z^{m+1}$ of~\eqref{eq:mnp-lambda} can be constructed in polynomial time from~\eqref{eq:lp}.
For every feasible basis~$B$ of~\eqref{eq:lp}, the following holds.
\begin{enumerate}
    \item[(a)] 
    The set~$\hat{B}=B\cup\{n+1\}$ is a corral.
    The improving atoms at~$\vec{w_{\hat{B}}}$ are those whose indices are improving at~$B$, with gaps~$g_j(\vec{w_{\hat{B}}})=\kappa\rc^B_j$ for~$j\in[n]$, for some~$\kappa>0$.
    \item[(b)]
    Let~$j\in\Gamma(B)$, and let~$B'=(B\setminus\{i^*\})\cup\{j\}$ be the pivot on~$j$.
    Then, the Wolfe major cycle that inserts~$\ivec{p}{j}$ at the corral~$\hat{B}$ only deletes~$i^*$, so the next corral is~$\widehat{B'}$.
    \item[(c)]
    The Wolfe runs from the corral~$\hat{B}$ are exactly the sequences~$\hat{B},\hat{B}_1,\dots,\hat{B}_f$ for the simplex runs~$B,B_1,\dots,B_f$ from~$B$.
    \item[(d)]
    Wolfe's method stops at~$\hat{B}$ with the optimal solution~$\vec{w_{\hat{B}}}$ of~\eqref{eq:mnp-lambda} if and only if~$B$ is optimal for~\eqref{eq:lp}.
\end{enumerate}
\end{theorem}

\ifsubmission\else
\section{Open problems}
\label{sec:open}
By Theorem~\ref{thm:equivalence}, an insertion rule that makes Wolfe's method polynomial from every corral would yield a polynomial simplex pivot rule for bounded nondegenerate LPs.
Conversely, it would be very interesting to prove a lower bound for all insertion rules simultaneously, as is known for the active-set method in convex quadratic maximization~\cite{BDHM26}.

\begin{openproblem}\label{op:rule}
Is there an insertion rule for which Wolfe's method needs only polynomially many major cycles in the number of atoms and the dimension?
\end{openproblem}

A lower bound for the linopt rule on base polytopes of submodular functions would yield the first lower bound for the Fujishige--Wolfe algorithm \cite{Fuj80,FI11}.

\begin{openproblem}\label{op:base}
Does Wolfe's method with the linopt rule need superpolynomially many major cycles on the base polytopes of some submodular functions?
\end{openproblem}
\fi

\ifsubmission\else
\subsubsection{\ackname}
The authors are grateful to L\'aszl\'o A.~V\'egh for pointing us to Wolfe's method as a potential application of our approach in~\cite{BDHM26}.
\fi

\begin{credits}
\subsubsection{AI Disclosure.}
This work was developed interactively with Claude Opus 5.5.
The starting point was to adapt the lower-bound constructions for the simplex method that are based on Markov decision processes~\cite{FHZ11,AF17,DFH23} to Wolfe's method.
The connection between Wolfe's method and the simplex method, on which Theorem~\ref{thm:equivalence} and all further results are based, emerged in discussions between the authors and the model.
The model further assisted in developing the proofs of all results, their formal verification in LEAN, and in the literature search. 
The paper was entirely rewritten by the authors, based on an initial draft by the model and using several LLMs for proof reading, typesetting, etc.
The authors verified the correctness of all statements and proofs and take full responsibility for the content of this paper.

\subsubsection{\discintname}
The authors have no competing interests to declare that are
relevant to the content of this article.
\end{credits}

\newpage

\bibliographystyle{splncs04}
\bibliography{references}

\clearpage
\appendix
\section{Useful facts in Section~\ref{sec:prelim}}
\label{app:prelim}

\subsection{The simplex method}
\label{app:simplex}

\begin{lemma}\label{lem:identities}
Let~$B$ be a basis of~\eqref{eq:lp}, and let~$j\in[n]$.
Then,
\begin{enumerate}
    \item[(a)] every~$i\in B$ has reduced cost~$\rc^B_i=0$,
    \item[(b)] the basic solution~$\ivec{x}{B}$ is the unique~$\vec{x}\in\Q^n$ with~$A\vec{x}=\vec{b}$ and~$x_k=0$ for all~$k\notin B$, and the objective value of~$B$ is~$z_B=\vec{b}^\top\ivec{y}{B}$,
    \item[(c)] if~$j$ is nonbasic, the basic direction~$\ivec{d}{B,j}$ is the unique~$\vec{d}\in\Q^n$ with~$A\vec{d}=\zero$,~$d_j=1$, and~$d_k=0$ for all~$k\notin B\cup\{j\}$, and its objective value is~$\vec{c}^\top\ivec{d}{B,j}=\rc^B_j$.
\end{enumerate}
\end{lemma}

\begin{proof}
\emph{Part (a).} By Definition~\ref{def:basis}, we have~$\rc^B_i=c_i-\vec{A_i}^\top\ivec{y}{B}=c_i-c_i=0$.

\emph{Part (b).} Let~$\vec{x}\in\R^n$ with~$x_k=0$ for all~$k\notin B$.
Then~$A\vec{x}=A_B\vec{x_B}$, so~$A\vec{x}=\vec{b}$ if and only if~$\vec{x_B}=A_B^{-1}\vec{b}$, since~$A_B$ is nonsingular.
The definition of the basic solution and Definition~\ref{def:basis_a} yield~$\vec{b}^\top\ivec{y}{B}=(\ivec{y}{B})^\top A_B\ivec{x}{B}_B=\vec{c_B}^\top\ivec{x}{B}_B=z_B$.

\emph{Part (c).} Let~$\vec{d}\in\R^n$ such that~$d_j=1$ and~$d_k=0$ for all~$k\notin B\cup\{j\}$.
Then, we have~$A\vec{d}=A_B\vec{d_B}+\vec{A_j}$, so~$A\vec{d}=\zero$ if and only if~$\vec{d_B}=-A_B^{-1}\vec{A_j}$, so~$\vec d$ is the basic direction.
With Definition~\ref{def:basis_a}, we get~$\vec{c}^\top\ivec{d}{B,j}=c_j-\vec{c_B}^\top A_B^{-1}\vec{A_j}=c_j-\vec{A_j}^\top\ivec{y}{B}=\rc^B_j$.
\qed
\end{proof}

\subsection{Wolfe's method}

\begin{proof}[of Lemma~\ref{lem:affine-min}]
\emph{Part (a).} This is Wolfe's criterion for the affine hull, we refer to~\cite[Lemma~1.2]{DHR} for a proof.

\emph{Part (b).} 
Let~$g_j(\vec{w_S})\neq 0$.
Assume~$\ivec{p}{j}\in\aff(S)$, then there are coefficients~$\vec \beta\in\R^S$ with~$\ivec{p}{j}=\sum_{i\in S}\beta_i\ivec{p}{i}$ and~$\sum_{i\in S}\beta_i=1$, so~$g_j(\vec{w_S})=\sum_{i\in S}\beta_ig_i(\vec{w_S})=0$ by part~(a).
This is a contradiction, so~$\ivec{p}{j}\notin\aff(S)$, and thus~$S\cup\{j\}$ is affinely independent.

\emph{Part (c).} Since~$S$ is a corral, we have~$\vec{w_S}\in\conv(S)$.
Then, the claim is Wolfe's criterion~\cite{Wolfe76}, see~\cite[Lemma~1.1]{DHR}.
\qed
\end{proof}

\section{Lemmas and proofs for Section~\ref{sec:equivalence}}
\label{app:proofs}

\subsection{Statements quoted in Section~\ref{sec:equivalence}}
\label{app:statements}

Every feasible solution~$\vec{\lambda}\in\R^{n+1}$ of~\eqref{eq:induced-lp} has objective value
\begin{equation}
    \vec{\gamma}^\top\vec{\lambda}\geq\one^\top\vec{\lambda}=1,
    \label{eq:positive-height}
\end{equation}
since~$\lambda_j\geq0$ and~$\gamma_j\geq1$, for every~$j\in[n+1]$, by~\eqref{eq:gamma}.

By~\eqref{eq:nf-dual}, every basis~$S$ of~\eqref{eq:induced-lp} satisfies $((\ivec{v}{S})^\top,h^S)M_S=\vec{\gamma_S}^\top$, which yields, since~$M$ has the columns $\vec{M_j}=((\ivec{\alpha}{j})^\top,1)^\top$, that
\begin{equation}
    (\ivec{v}{S})^\top\ivec{\alpha}{i}+h^S=\gamma_i\qquad\text{for every }i\in S.
    \label{eq:dual}
\end{equation}

\begin{lemma}\label{lem:T-bound}
Consider an admissible linear program, and let~$U\geq1$ be the largest absolute value of an entry of~$A$ or~$\vec{b}$.
For every feasible point~$\vec{x}\in\polytope$, we have $2\cdot\one^\top\vec{x}<T$, where~$T=2m(mU)^m+1$.
\end{lemma}
\begin{proof}
Let~$B$ be a feasible basis.
By Cramer's rule, every basic coordinate is $\icoord{x}{B}{i}=\det A_B^{(i)}/\det A_B$, where~$A_B^{(i)}$ arises from~$A_B$ by replacing the column of~$i$ by~$\vec{b}$.
Both matrices are integral of order~$m$, and their entries have absolute value at most~$U$.
Hence~$\det A_B$ is a nonzero integer, so~$|\det A_B|\geq1$.
Every column of~$A_B^{(i)}$ has norm at most~$\sqrt mU$, so Hadamard's inequality, which bounds the absolute value of a determinant by the product of the norms of its columns, gives~$|\det A_B^{(i)}|\leq(\sqrt mU)^m\leq(mU)^m$.
Since~$\ivec{x}{B}\geq\zero$ is zero outside the~$m$ indices of~$B$, we get $2\cdot\one^\top\ivec{x}{B}\leq2m(mU)^m<T$.
By Section~\ref{sec:simplex}, the maximum of $\one^\top\vec{x}$ over~$\polytope$ is attained at such a basic solution $\ivec{x}{B}$.
\qed
\end{proof}

\begin{lemma}\label{lem:induced-facts}
Let~$S$ be a basis of~\eqref{eq:induced-lp}, and let~$j\in[n+1]$.
Then,
\begin{enumerate}
    \item[(a)] we have $\sum_{i\in S}\icoord{\lambda}{S}{i}\ivec{\alpha}{i}=\zero$ and $\sum_{i\in S}\icoord{\lambda}{S}{i}=1$,
    \item[(b)] the objective value of~$S$ is $\vec{\gamma}^\top\ivec{\lambda}{S}=h^S$; if~$S$ is feasible, then~$h^S\geq1$,
    \item[(c)] if~$j$ is nonbasic, $\sum_{k\in[n+1]}\icoord{d}{S,j}{k}\ivec{\alpha}{k}=\zero$, $\sum_{k\in[n+1]}\icoord{d}{S,j}{k}=0$, and $\vec{\gamma}^\top\ivec{d}{S,j}=\rg^S_j$.
\end{enumerate}
\end{lemma}
\begin{proof}
\emph{Part~(a).}
By~\eqref{eq:nf-solution} and as~$M$ has the columns~$\vec{M_i}=((\ivec{\alpha}{i})^\top,1)^\top$, we have
\[
    (\zero^\top,1)^\top=M_S\ivec{\lambda}{S}_S=\sum_{i\in S}\icoord{\lambda}{S}{i}\vec{M_i}=\Bigl(\sum_{i\in S}\icoord{\lambda}{S}{i}(\ivec{\alpha}{i})^\top,\ \sum_{i\in S}\icoord{\lambda}{S}{i}\Bigr)^\top.
\]

\emph{Part~(b).}
By~\eqref{eq:nf-solution} and~\eqref{eq:nf-dual}, we have
\[
    \vec{\gamma}^\top\ivec{\lambda}{S}=\vec{\gamma_S}^\top M_S^{-1}(\zero^\top,1)^\top=((\ivec{v}{S})^\top,h^S)(\zero^\top,1)^\top=h^S.
\]
If~$S$ is feasible, then~\eqref{eq:positive-height} gives~$h^S\geq1$.

\emph{Part~(c).}
Write $\vec{d}=\ivec{d}{S,j}$.
By Definition~\ref{def:basis_d}, applied to~\eqref{eq:induced-lp}, we have~$d_j=1$ and~$\vec{d_S}=-M_S^{-1}\vec{M_j}$, and all other coordinates of~$\vec{d}$ are zero.
Hence
\[
    \Bigl(\sum_{k\in[n+1]}d_k(\ivec{\alpha}{k})^\top,\ \sum_{k\in[n+1]}d_k\Bigr)^\top
    =\sum_{k\in[n+1]}d_k\vec{M_k}
    =M_S\vec{d_S}+\vec{M_j}
    =\zero.
\]
We use~\eqref{eq:nf-direction}, \eqref{eq:nf-dual}, and~\eqref{eq:nf-rc} to compute the objective value
\begin{align*}
    \vec{\gamma}^\top\vec{d}&{}=\gamma_j+\vec{\gamma_S}^\top\vec{d_S}=\gamma_j-\vec{\gamma_S}^\top M_S^{-1}\vec{M_j}\\
    &=\gamma_j-((\ivec{v}{S})^\top,h^S)\vec{M_j}=\gamma_j-(\ivec{v}{S})^\top\ivec{\alpha}{j}-h^S=\rg^S_j.\tag*{$\qed$}
\end{align*}
\end{proof}

\begin{lemma}\label{lem:encoding}
Let~$(A,\vec{b},\vec{c})$ be the data of the admissible linear program~\eqref{eq:lp}.
Then, the atoms $\ivec{p}{j}$ defined in~\eqref{eq:atoms} are integer vectors.
The numbers~$T$,~$C$,~$K$,~$H$,~$N$ defined in~\eqref{eq:T},~\eqref{eq:gamma}, and~\eqref{eq:scale}
and the atoms $\ivec{p}{j}$ have bit lengths that are polynomial in the bit length of~$(A,\vec{b},\vec{c})$, and they can be computed in polynomial time.
\end{lemma}
\begin{proof}
By~\ref{ass:integral}, \eqref{eq:alpha}, and~\eqref{eq:gamma}, the objective coefficients~$\gamma_{n+1}=C$ and~$\gamma_j=C+c_j$ for~$j\in[n]$ and the vectors $T\ivec{\alpha}{j}$ for~$j\in[n+1]$ are integral.
By~\eqref{eq:scale}, the number~$N$ is an integer multiple of~$T$.
Hence, every atom
\[
    \ivec{p}{j}=\bigl(N\gamma_j,\ (N^2/T)\,T(\ivec{\alpha}{j})^\top\bigr)^\top,\qquad\text{for }j\in[n+1],
\]
from~\eqref{eq:atoms} is an integer vector.

Let~$U\geq1$ be the largest absolute value of an entry of~$A$ or~$\vec{b}$, as in Lemma~\ref{lem:T-bound}.
By~\eqref{eq:T}, we have~$\log T=O(m\log(mU))$.
The number~$K$ is the largest absolute value of the entries of the integral vectors $(\gamma_j,T(\ivec{\alpha}{j})^\top,1)^\top$ for~$j\in[n+1]$.
We have~$|\gamma_j|\leq C+\max_{k\in[n]}|c_k|$ by~\eqref{eq:gamma}, and by~\eqref{eq:alpha}, the entries of $T\ivec{\alpha}{j}=T\vec{A_j}-\vec{b}$ for~$j\in[n]$ and of $T\ivec{\alpha}{n+1}=-\vec{b}$ have absolute value at most~$(T+1)U$.
Hence
\[
    K\leq\max\Bigl\{C+\max_{k\in[n]}|c_k|,\ (T+1)U\Bigr\},
\]
and~$\log K$ is polynomial in the bit length of~$(A,\vec{b},\vec{c})$, since~$C=1+\max_{k\in[n]}|c_k|$.
By~\eqref{eq:scale}, we have
\[
    \log H=(m+1)\log\bigl((m+1)K^2\bigr)\quad\text{and}\quad\log N=\log(2mT)+3\log H,
\]
so~$\log H$ and~$\log N$ are polynomial as well.
Every coordinate of an atom has absolute value at most~$N^2K$, and there are~$(m+1)(n+1)$ such coordinates.
Thus, the atoms also have polynomial bit lengths.

Finally, evaluating the formulas~\eqref{eq:T}, \eqref{eq:gamma}, \eqref{eq:scale}, and~\eqref{eq:atoms} for~$T$,~$C$,~$K$,~$H$,~$N$, and the atoms takes~$O(mn)$ additions, multiplications, and comparisons, and every intermediate result is an integer of absolute value at most~$N^2K$.
\qed
\end{proof}

\subsection{Proofs for Section~\ref{sec:equivalence}}
\label{app:proofs-equivalence}

\begin{proof}[of Lemma~\ref{lem:hom}]
Let~$\bar{A}\in\Q^{(m+1)\times(n+1)}$ be the matrix with columns~$\vec{\bar{A}_j}=(\vec{A_j}^\top,1)^\top$ for~$j\in[n]$ and~$\vec{\bar{A}_{n+1}}=(\zero^\top,1)^\top$, which is the constraint matrix of~\eqref{eq:lp} with the added constraint~$\one^\top\vec{x}+x_{n+1}=T$.
By~\eqref{eq:alpha}, we have
\begin{equation}
    M=E\bar{A}\quad\text{and}\quad(\zero^\top,1)^\top=E\,(\vec{b}^\top/T,1)^\top,\qquad\text{where}\quad
    E=\begin{pmatrix}I_m&-\frac{1}{T}\vec{b}\\ \zero^\top&1\end{pmatrix}.
    \label{eq:factor}
\end{equation}
Since~$E$ is nonsingular, every~$\vec{\lambda}\in\R^{n+1}$ and every~$\vec{d}\in\R^{n+1}$ satisfy
\begin{equation}
    \begin{aligned}
        M\vec{\lambda}=(\zero^\top,1)^\top&\iff A\bigl(T\vec{\lambda_{[n]}}\bigr)=\vec{b}\ \text{ and }\ \lambda_{n+1}=1-\one^\top\vec{\lambda_{[n]}},\\
        M\vec{d}=\zero&\iff A\vec{d_{[n]}}=\zero\ \text{ and }\ d_{n+1}=-\one^\top\vec{d_{[n]}}.
    \end{aligned}
    \label{eq:hom-systems}
\end{equation}
The lemma follows from these two equivalences by standard textbook arguments, which we include for completeness.

Let~$B\subseteq[n]$ have~$m$ elements.
The column of~$n+1$ in~$\bar{A}$ is the last unit vector, so expanding along it gives~$\det\bar{A}_{\hat{B}}=\pm\det A_B$, and~$\det M_{\hat{B}}=\det\bar{A}_{\hat{B}}$ by~\eqref{eq:factor}.
Hence~$\hat{B}$ is a basis of~\eqref{eq:induced-lp} if and only if~$B$ is a basis of~\eqref{eq:lp}; in particular,~$M$ has full row rank.
Let~$B$ be a basis of~\eqref{eq:lp}.
The vector~$\vec{\lambda}\coloneqq\bigl((\ivec{x}{B})^\top,\,T-\one^\top\ivec{x}{B}\bigr)^\top/T$ vanishes outside~$\hat{B}$ and satisfies~$M\vec{\lambda}=(\zero^\top,1)^\top$ by~\eqref{eq:hom-systems}, so~$\ivec{\lambda}{\hat{B}}=\vec{\lambda}$ by Lemma~\ref{lem:identities}(b), applied to~\eqref{eq:induced-lp}.
This proves part~(a) of the statement.
If~$B$ is feasible, then all coordinates of~$\ivec{\lambda}{\hat{B}}$ in~$\hat{B}$ are positive by~\ref{ass:nondeg} and Lemma~\ref{lem:T-bound}, so~$\hat{B}$ is feasible and nondegenerate.
Conversely, let~$S$ be a feasible basis of~\eqref{eq:induced-lp}, and let~$\vec{x}\coloneqq T\ivec{\lambda}{S}_{[n]}$.
Then~$\vec{x}\in\polytope$ by~\eqref{eq:hom-systems}, so~$\icoord{\lambda}{S}{n+1}=1-\one^\top\vec{x}/T>0$ by Lemma~\ref{lem:T-bound}.
Hence~$n+1\in S$, and~$B\coloneqq S\setminus\{n+1\}$ is a basis of~\eqref{eq:lp} with~$\hat{B}=S$.
Since~$\vec{x}$ vanishes outside~$B$, Lemma~\ref{lem:identities}(b) gives~$\ivec{x}{B}=\vec{x}\geq\zero$, so~$B$ is feasible.
Thus~$B\mapsto\hat{B}$ is a bijection between the feasible bases, and all feasible bases of~\eqref{eq:induced-lp} are nondegenerate.
Moreover, the feasible set of~\eqref{eq:induced-lp} is nonempty, since $\hat{B}$ is feasible for every feasible basis $B$ of~\eqref{eq:lp}, and such a basis exists by~\ref{ass:bounded}.
It is also bounded, since the entries of every feasible point are nonnegative and sum to one.
This proves the first sentence of the lemma.

We now prove parts~(b) and~(c).
Let~$B$ be a basis of~\eqref{eq:lp}, let~$\vec{v}\coloneqq\ivec{v}{\hat{B}}$ and~$h\coloneqq h^{\hat{B}}$, and let~$(\vec{y'}^\top,h')\coloneqq(\vec{v}^\top,h)E=(\vec{v}^\top,\,h-\vec{v}^\top\vec{b}/T)$.
By definition of the dual vector~\eqref{eq:nf-dual} and~\eqref{eq:factor}, we have~$(\vec{y'}^\top,h')\bar{A}_{\hat{B}}=(\vec{v}^\top,h)M_{\hat{B}}=\vec{\gamma_{\hat{B}}}^\top$.
The column of~$n+1$ gives~$h'=\gamma_{n+1}=C$, so the columns of the indices~$i\in B$ give~$\vec{A_i}^\top\vec{y'}=\gamma_i-C=c_i$.
This gives~$\vec{y'}=\ivec{y}{B}$ by Definition~\ref{def:basis_a}.
Hence~$\vec{v}=\ivec{y}{B}$, and Lemma~\ref{lem:identities}(b) gives~$h=C+\vec{b}^\top\ivec{y}{B}/T=C+z_B/T$, which proves part~(b).
By~\eqref{eq:nf-rc}, every index~$j\in[n]$ satisfies
\[
    \rg^{\hat{B}}_j=\gamma_j-(\vec{v}^\top,h)\vec{M_j}=\gamma_j-(\vec{y'}^\top,h')\vec{\bar{A}_j}=C+c_j-\vec{A_j}^\top\ivec{y}{B}-C=\rc^B_j,
\]
and~$\rg^{\hat{B}}_{n+1}=0$ by Lemma~\ref{lem:identities}(a), applied to~\eqref{eq:induced-lp}.
This proves part~(c).

Let~$B$ be a feasible basis of~\eqref{eq:lp}, let~$j\in\Gamma(B)$, and let~$\vec{d}\coloneqq\ivec{d}{B,j}$.
The vector~$(\vec{d}^\top,-\one^\top\vec{d})^\top$ has the coordinate~$1$ at~$j$, vanishes outside~$\hat{B}\cup\{j\}$, and lies in the kernel of~$M$ by~\eqref{eq:hom-systems} and Lemma~\ref{lem:identities}(c).
So it is~$\ivec{d}{\hat{B},j}$ by Lemma~\ref{lem:identities}(c), applied to~\eqref{eq:induced-lp}.
Together with part~(a), this shows that the ratio of every index~$i\in B$ in the ratio test~\eqref{eq:ratio-induced} is its ratio in~\eqref{eq:ratio} divided by~$T$.
Since~$j\notin\hat{B}$ and~$\rg^{\hat{B}}_j=\rc^B_j<0$ by part~(c), the statements after Definition~\ref{def:basis}, applied to~\eqref{eq:induced-lp}, show that the pivot on~$j$ at~$\hat{B}$ has a unique leaving index and gives a feasible basis, which contains~$n+1$ by the bijection above.
So the leaving index lies in~$B$ and attains the minimum in~\eqref{eq:ratio}, and it is therefore the leaving index~$i^*$ of the pivot on~$j$ at~$B$.

By part~(c), the improving indices at~$\hat{B}$ are those at~$B$, so the simplex method stops at~$\hat{B}$ if and only if it stops at~$B$, and the pivots on~$j$ at~$B$ and at~$\hat{B}$ give~$B'$ and~$\widehat{B'}$, with the same leaving index.
Since~$B\mapsto\hat{B}$ is a bijection between the feasible bases, induction on the number of pivots shows that the simplex runs correspond as stated.
\qed
\end{proof}

\begin{proof}[of Lemma~\ref{lem:direction}]
Write $\vec{d}=\ivec{d}{S,j}$.
By Definition~\ref{def:basis_d}, applied to the program~\eqref{eq:induced-lp}, we have~$d_j=1$ and~$d_k=0$ for every~$k\in[n+1]\setminus(S\cup\{j\})$.
Hence, Lemma~\ref{lem:induced-facts}(c) gives
\[
    \ivec{\alpha}{j}+\sum_{i\in S}d_i\ivec{\alpha}{i}=\zero,\qquad
    1+\sum_{i\in S}d_i=0,\qquad
    \gamma_j+\sum_{i\in S}d_i\gamma_i=\rg^S_j.
\]
By the definition~\eqref{eq:atoms} of the atoms, the first and the third equation yield
\[
    \ivec{p}{j}+\sum_{i\in S}d_i\ivec{p}{i}
    =\Bigl(N\Bigl(\gamma_j+\sum_{i\in S}d_i\gamma_i\Bigr),\ N^2\Bigl(\ivec{\alpha}{j}+\sum_{i\in S}d_i\ivec{\alpha}{i}\Bigr)^\top\Bigr)^\top
    =(N\rg^S_j,\zero^\top)^\top
    =N\rg^S_j\vec{e_1},
\]
which is the first identity.
Let~$\vec{w}\in\R^{m+1}$.
We take the inner product of the first identity with~$\vec{w}$ and subtract~$\bigl(1+\sum_{i\in S}d_i\bigr)\norm{\vec{w}}^2=0$.
By the definition of the gaps (Definition~\ref{def:gap}), this yields
\[
    g_j(\vec{w})+\sum_{i\in S}d_ig_i(\vec{w})
    =\vec{w}^\top\ivec{p}{j}-\norm{\vec{w}}^2+\sum_{i\in S}d_i\bigl(\vec{w}^\top\ivec{p}{i}-\norm{\vec{w}}^2\bigr)
    =N\rg^S_j\vec{e_1}^\top\vec{w}
    =N\rg^S_jw_1.\tag*{$\qed$}
\]
\end{proof}

\begin{proof}[of Lemma~\ref{lem:corral}]
\emph{Part~(a).}
Let~$\vec{\beta}\in\R^S$ with $\sum_{i\in S}\beta_i\ivec{p}{i}=\zero$ and~$\sum_{i\in S}\beta_i=0$.
By the definition of the atoms~\eqref{eq:atoms}, the constraint coordinates give $N^2\sum_{i\in S}\beta_i\ivec{\alpha}{i}=\zero$, so~$M_S\vec{\beta}=\zero$.
Since~$S$ is a basis,~$M_S$ is nonsingular, so~$\vec{\beta}=\zero$.
Hence, the atoms~$\ivec{p}{i}$ for~$i\in S$ are affinely independent.
Write~$\varepsilon=\varepsilon^S$,~$s=s^S$, $\vec{v}=\ivec{v}{S}$,~$h=h^S$, and let~$\vec{\lambda}\coloneqq\ivec{\lambda}{S}-\varepsilon\ivec{u}{S}$.
Then $M_S\vec{\lambda_S}=(\zero^\top,1)^\top-\varepsilon(\vec{v}^\top,0)^\top=(-\varepsilon\vec{v}^\top,1)^\top$, so the coordinates of~$\vec{\lambda}$ sum to~$1$ and $\sum_{i\in S}\lambda_i\ivec{\alpha}{i}=-\varepsilon\vec{v}$.
By~\eqref{eq:projection-quantities}, we have
\[
    N^2\varepsilon
    =\frac{N^2h}{N^2+\norm{\vec{v}}^2}
    =h-\frac{h}{N^2+\norm{\vec{v}}^2}\norm{\vec{v}}^2
    =h-\varepsilon\norm{\vec{v}}^2.
\]
With Lemma~\ref{lem:induced-facts}(b) and~\eqref{eq:nf-dual}, which gives $\vec{\gamma_S}^\top M_S^{-1}=(\vec{v}^\top,h)$, we obtain
\[
    \vec{\gamma}^\top\vec{\lambda}
    =\vec{\gamma_S}^\top\ivec{\lambda}{S}_S-\varepsilon\vec{\gamma_S}^\top M_S^{-1}(\vec{v}^\top,0)^\top
    =h-\varepsilon(\vec{v}^\top,h)(\vec{v}^\top,0)^\top
    =h-\varepsilon\norm{\vec{v}}^2
    =N^2\varepsilon
    =s.
\]
We define $\vec{w}\coloneqq\sum_{i\in S}\lambda_i\ivec{p}{i}$.
By the penalty identity~\eqref{eq:penalty}, we have $\vec{w}=(Ns,-N^2\varepsilon\vec{v}^\top)^\top$.
For~$i\in S$, the dual equation~\eqref{eq:dual} and~$s=N^2\varepsilon=h-\varepsilon\norm{\vec{v}}^2$ give
\begin{align*}
    \vec{w}^\top\ivec{p}{i}&=N^2s\gamma_i-N^4\varepsilon\vec{v}^\top\ivec{\alpha}{i}=N^4\varepsilon(\gamma_i-\vec{v}^\top\ivec{\alpha}{i})=N^4\varepsilon h,\\
    \norm{\vec{w}}^2&=N^2s^2+N^4\varepsilon^2\norm{\vec{v}}^2=N^4\varepsilon\bigl(s+\varepsilon\norm{\vec{v}}^2\bigr)=N^4\varepsilon h.
\end{align*}
So~$g_i(\vec{w})=0$ for every~$i\in S$, and Lemma~\ref{lem:affine-min}(a) gives~$\vec{w_S}=\vec{w}$ and $\ivec{\bar{\lambda}}{S}=\vec{\lambda}$.

\emph{Part~(b).}
By Lemma~\ref{lem:affine-min}(a), we have~$g_i(\vec{w_S})=0$ for every~$i\in S$, and the objective coordinate of~$\vec{w_S}$ is~$Ns^S$ by part~(a).
For every index~$j\in[n+1]\setminus S$, Lemma~\ref{lem:direction} gives~$g_j(\vec{w_S})=N\rg^S_j\cdot Ns^S$; for every index~$j\in S$, both sides vanish, since~$\rg^S_j=0$ by Lemma~\ref{lem:identities}(a) for~\eqref{eq:induced-lp}.
If~$S$ is feasible, then~$h^S\geq1$ by Lemma~\ref{lem:induced-facts}(b), and hence we obtain~$s^S=N^2h^S/(N^2+\norm{\ivec{v}{S}}^2)>0$.
\qed
\end{proof}

\begin{proof}[of Lemma~\ref{lem:bases-corrals}]
The entries of~$M$ are rational, so we first make them integral, as Cramer's rule and Hadamard's inequality need.
Let~$\widetilde{M}$ arise from~$M$ by multiplying its first~$m$ rows by~$T$.
Its columns are the integer vectors $(T(\ivec{\alpha}{k})^\top,1)^\top$ for~$k\in[n+1]$, whose entries have absolute value at most~$K$ by the definition of~$K$ before~\eqref{eq:scale}.
Since~$S$ is a basis,~$\widetilde{M}_S$ is nonsingular, and the basic solution, the basic directions, the dual vector, and $\ivec{u}{S}$ from~\eqref{eq:projection-quantities}
solve the systems
\begin{equation}
    \begin{aligned}
        \widetilde{M}_S\ivec{\lambda}{S}_S&=(\zero^\top,1)^\top,&\qquad
        \widetilde{M}_S\ivec{d}{S,j}_S&=-(T(\ivec{\alpha}{j})^\top,1)^\top,\\
        \widetilde{M}_S^\top\bigl((\ivec{v}{S})^\top/T,\,h^S\bigr)^\top&=\vec{\gamma_S},&
        \widetilde{M}_S\ivec{u}{S}&=(T(\ivec{v}{S})^\top,0)^\top.
    \end{aligned}
    \label{eq:integer-systems}
\end{equation}
Let~$\widetilde{M}_S\vec{z}=\vec{r}$ or $\widetilde{M}_S^\top\vec{z}=\vec{r}$ be a system with an integer right-hand side~$\vec{r}$ with~$\norm{\vec{r}}_\infty\leq K$.
By Cramer's rule, every entry of~$\vec{z}$ is a quotient of two determinants of integer matrices of order~$m+1$ with entries of absolute value at most~$K$, and the denominator~$\det\widetilde{M}_S$ is nonzero.
By Hadamard's inequality, these determinants have absolute value at most~$(\sqrt{m+1}K)^{m+1}\leq H$, so
\begin{equation}
    |z_\ell|\leq H\quad\text{for every }\ell,\qquad\text{and}\quad|z_\ell|\geq1/H\quad\text{if }z_\ell\neq0.
    \label{eq:det-bound}
\end{equation}
This applies to the first three systems in~\eqref{eq:integer-systems}, whose right-hand sides are integer vectors with entries of absolute value at most~$K$ by the definition of~$K$, and it applies to the systems~$\widetilde{M}_S\vec{z}=\vec{e_\ell}$, whose solutions are the columns of~$\widetilde{M}_S^{-1}$.

\emph{Part~(a).}
The basic coordinates of $\ivec{\lambda}{S}$ are nonzero, because~$S$ is non\-degenerate by Lemma~\ref{lem:hom}, so they are at least $\icoord{\lambda}{S}{i} \geq 1/H$ by~\eqref{eq:det-bound}.

\emph{Part~(b).}
We have $\icoord{d}{S,j}{j}=1$ by Definition~\ref{def:basis_d}, the coordinates of $\ivec{d}{S,j}$ in~$S$ satisfy~\eqref{eq:det-bound} by the second system in~\eqref{eq:integer-systems}, and all other coordinates are zero.

\emph{Part~(c).}
The equality is Lemma~\ref{lem:corral}(a), since~$\varepsilon^S>0$.
By~\eqref{eq:det-bound} for the third system in~\eqref{eq:integer-systems}, we have $\norm{\ivec{v}{S}}_\infty\leq TH$ and~$h^S\leq H$.
By~\eqref{eq:det-bound} for the columns of~$\widetilde{M}_S^{-1}$, every entry of~$\widetilde{M}_S^{-1}$ has absolute value at most~$H$.
Since $\ivec{u}{S}=\widetilde{M}_S^{-1}(T(\ivec{v}{S})^\top,0)^\top$ by the fourth system, we obtain $\norm{\ivec{u}{S}}_\infty\leq m\cdot H\cdot T^2H{}=mT^2H^2$.
Since~$\varepsilon^S\leq h^S/N^2\leq H/N^2$ by~\eqref{eq:projection-quantities}, the choice~\eqref{eq:scale} of~$N$ gives
\[
    \varepsilon^S\norm{\ivec{u}{S}}_\infty
       \leq\frac{mT^2H^3}{N^2}
       =\frac{mT^2H^3}{4m^2T^2H^6}
       =\frac{1}{4mH^3}.
\]
Finally, (a) and (c) give $\icoord{\bar{\lambda}}{S}{i}\geq1/H-1/(4mH^3)>0$ for every~$i\in S$.
\qed
\end{proof}

\begin{proof}[of Lemma~\ref{lem:pivot}]
Write $\vec{d}=\ivec{d}{S,j}$,~$\theta=\theta^{S,j}$, and $\rho=\varepsilon^S\norm{\ivec{u}{S}}_\infty$, and extend all vectors on~$S$ and on~$S'$ by zero to~$S\cup\{j\}$; in particular, $\icoord{u}{S}{j}=0$.
Recall that~$d_j=1$ by Definition~\ref{def:basis_d} and that~$d_{i^*}<0$, because the minimum in the ratio test~\eqref{eq:ratio-induced} is taken over the indices~$i$ with~$d_i<0$.
By~\eqref{eq:pivot}, which applies to~\eqref{eq:induced-lp} by Lemma~\ref{lem:hom}, and by Lemma~\ref{lem:induced-facts}(b), the set~$S'$ is a feasible basis with
\[
    \ivec{\lambda}{S'}=\ivec{\lambda}{S}+\theta\vec{d}\quad\text{and}\quad h^{S'}=h^S+\theta\rg^S_j.
\]
By Lemma~\ref{lem:bases-corrals}, both~$S$ and~$S'$ are corrals.
We follow the lines of Algorithm~\ref{alg:wolfe}.

\emph{Insertion (Lines~\ref{line:wf-insert} and~\ref{line:wf-loop}).}
Line~\ref{line:wf-insert} inserts~$j$ with coefficient~$0$, so the current coefficient vector is $\ivec{\bar{\lambda}}{S}$.
By Lemma~\ref{lem:corral}(b), we have~$g_j(\vec{w_S})=N^2s^S\rg^S_j<0$, so the~$m+2$ atoms of~$S\cup\{j\}$ are affinely independent by Lemma~\ref{lem:affine-min}(b).
Hence, their affine hull is~$\R^{m+1}$, and their affine minimizer is the origin.
Let
\[
    \tau\coloneqq\frac{h^S}{-\rg^S_j}>0\quad\text{and}\quad\vec{\mu}\coloneqq\ivec{\lambda}{S}+\tau\vec{d}.
\]
By Lemma~\ref{lem:induced-facts}(a)--(c), we have $M\vec{\mu}=(\zero^\top,1)^\top$ and $\vec{\gamma}^\top\vec{\mu}=h^S+\tau\rg^S_j=0$.
So the coordinates of~$\vec{\mu}$ sum to~$1$, and the penalty identity~\eqref{eq:penalty} gives $\sum_{k\in S\cup\{j\}}\mu_k\ivec{p}{k}=\zero$.
Since affine coefficients with respect to affinely independent atoms are unique, $\vec{\mu}=\ivec{\bar{\lambda}}{S\cup\{j\}}$ is the coefficient vector of the affine minimizer of~$S\cup\{j\}$.
Since~$\rg^S_j<0$, the objective value~$h^S+r\rg^S_j$ strictly decreases in~$r$; it is~$h^{S'}\geq1$ at~$r=\theta$ by Lemma~\ref{lem:induced-facts}(b) and~$0$ at~$r=\tau$, so~$\theta<\tau$.
With $\icoord{\lambda}{S'}{i^*}=0$, this gives
\[
    \mu_{i^*}=\icoord{\lambda}{S}{i^*}+\tau d_{i^*}=\icoord{\lambda}{S'}{i^*}+(\tau-\theta)d_{i^*}=(\tau-\theta)d_{i^*}<0.
\]
So~$S\cup\{j\}$ is not a corral, and the loop in Line~\ref{line:wf-loop} performs a minor cycle.

\emph{Minor cycle (Lines~\ref{line:wf-direction}--\ref{line:wf-delete}).}
Line~\ref{line:wf-direction} sets the direction $\vec{\mu}-\ivec{\bar{\lambda}}{S}$, so by Lemma~\ref{lem:corral}(a), Lines~\ref{line:wf-step} and~\ref{line:wf-move} move along the segment
\[
    \vec{\eta(t)}\coloneqq(1-t)\ivec{\bar{\lambda}}{S}+t\vec{\mu}=\vec{\xi(t)}-(1-t)\varepsilon^S\ivec{u}{S}\qquad\text{for }t\in[0,1]
\]
up to the largest~$t$ with~$\vec{\eta(t)}\geq\zero$, where $\vec{\xi(t)}\coloneqq\ivec{\lambda}{S}+t\tau\vec{d}$.
The point~$\vec{\xi(t)}$ lies on the simplex ray, every coordinate of~$\vec{\eta(t)}$ differs from that of~$\vec{\xi(t)}$ by at most~$\rho$, and the ray reaches $\ivec{\lambda}{S'}$ at~$t_0\coloneqq\theta/\tau\in(0,1)$.
The coordinate~$\eta_{i^*}$ is an affine function of~$t$, which is positive at~$0$, since~$S$ is a corral, and negative at~$1$, since~$\mu_{i^*}<0$.
Let~$\hat t\in(0,1)$ be its zero.
Since~$\ivec{\lambda}{S'}=\ivec{\lambda}{S}+\theta\vec{d}$ and~$\theta=t_0\tau$, every~$t\in[0,1]$ satisfies
\[
    \vec{\xi(t)}
    =\ivec{\lambda}{S}+t\tau\vec{d}
    =\ivec{\lambda}{S}+\theta\vec{d}+(t\tau-\theta)\vec{d}
    =\ivec{\lambda}{S'}+(t-t_0)\tau\vec{d}.
\]
Since~$\icoord{\lambda}{S'}{i^*}=0$, this gives~$\xi_{i^*}(\hat t)=(\hat t-t_0)\tau d_{i^*}$.
Therefore, since~$\eta_{i^*}(\hat t)=0$, the bound~$|\xi_{i^*}(\hat t)|=|\xi_{i^*}(\hat t)-\eta_{i^*}(\hat t)|\leq\rho$ now gives~$|\hat t-t_0|\tau\leq\rho/|d_{i^*}|$.

Now let~$i\in S\setminus\{i^*\}$.
By Lemma~\ref{lem:bases-corrals}(a) for~$S'$, we have $\xi_i(t_0)=\icoord{\lambda}{S'}{i}\geq1/H$.
By Lemma~\ref{lem:bases-corrals}(b) and~(c), we have~$|d_i|\leq H^2|d_{i^*}|$ and~$\rho\leq1/(4mH^3)$.
Since~$1+H^2\leq2H^2$, we obtain
\[
    \eta_i(\hat t)\geq\xi_i(\hat t)-\rho=\icoord{\lambda}{S'}{i}+(\hat t-t_0)\tau d_i-\rho
    \geq\frac1H-\rho\Bigl(\frac{|d_i|}{|d_{i^*}|}+1\Bigr)
    \geq\frac1H-\frac{2H^2}{4mH^3}>0.
\]
Since~$S$ is a corral, also $\eta_i(0)=\icoord{\bar{\lambda}}{S}{i}>0$, so the affine function~$\eta_i$ is positive on~$[0,\hat t]$.
Moreover, the coordinate~$\eta_j(t)=t\tau$ is positive for~$t>0$, and~$\eta_{i^*}$ is nonnegative on~$[0,\hat t]$ and negative after~$\hat t$.
Hence Line~\ref{line:wf-step} gives the step length~$\hat t$, and the coordinate of~$i^*$ is the only zero coordinate of~$\vec{\eta(\hat t)}$ in~$S\cup\{j\}$.
So Line~\ref{line:wf-delete} deletes exactly~$i^*$, and the new support is~$S'$.

\emph{Acceptance (Lines~\ref{line:wf-loop} and~\ref{line:wf-accept}).}
Since~$S'$ is a corral, the loop in Line~\ref{line:wf-loop} ends, Line~\ref{line:wf-accept} moves to the coefficients of the affine minimizer of~$S'$, and the major cycle ends with the corral~$S'$.
\qed
\end{proof}

\begin{proof}[of Theorem~\ref{thm:equivalence}]
The guarantee on the construction time and the integrality of the atoms are given by Lemma~\ref{lem:encoding}.

\emph{Part~(a).} By Lemma~\ref{lem:hom}, the set~$\hat{B}$ is a feasible basis of~\eqref{eq:induced-lp}, so it is a corral by Lemma~\ref{lem:bases-corrals}.
The gaps follow from Lemma~\ref{lem:corral}(b) and Lemma~\ref{lem:hom}, with~$\kappa=N^2s^{\hat{B}}>0$.

\emph{Part~(b).} By Lemma~\ref{lem:hom}, we have~$\rg^{\hat{B}}_j=\rc^B_j<0$ and the pivot on~$j$ at~$\hat{B}$ has the leaving index~$i^*$ and hence gives~$\widehat{B'}$.
Thus, Lemma~\ref{lem:pivot} gives the claim.

\emph{Part~(c).} 
By~(a), we have~$\Gamma(\vec{w_{\hat{B}}})=\Gamma(B)$, since~$n+1$ has gap zero by Lemma~\ref{lem:affine-min}(a).
So Wolfe's method stops at~$\hat{B}$ if and only if the simplex method stops at~$B$, and otherwise, it may insert exactly the indices~$j\in\Gamma(B)$ that may enter at~$B$.
By~(b), the major cycle that inserts such an index~$j$ ends at the corral~$\widehat{B'}$ for the pivot~$B'$ on~$j$ at~$B$, which is again a feasible basis.
Let~$B_0\coloneqq B$.
Induction on~$k$ shows that a sequence~$S_0,\dots,S_f$ with~$S_0=\hat{B}_0$ is a Wolfe run if and only if there is a simplex run~$B_0,\dots,B_f$ from~$B_0$ with~$S_k=\hat{B}_k$ for every~$k\in\{0,\dots,f\}$.

\emph{Part~(d).} By~(a), no atom is improving at~$\vec{w_{\hat{B}}}$ if and only if~$\Gamma(B)=\emptyset$, in which case~$B$ is optimal.
Further, Lemma~\ref{lem:affine-min}(c) applies in this case, since~$\hat{B}$ is a corral by~(a), so~$\vec{w_{\hat{B}}}$ is optimal for~\eqref{eq:mnp-lambda}.
\qed
\end{proof}

\section{Corollaries of Theorem \ref{thm:equivalence}}
\label{app:consequences}

\subsection{Proof of Corollary~\ref{cor:intro-polynomial}}
\label{app:cor-polynomial}

\begin{proof}[of Corollary~\ref{cor:intro-polynomial}]
Consider an insertion rule that needs polynomially many major cycles from every starting corral, and an admissible linear program with a feasible basis~$B_0$.
A pivot rule can construct the atoms of Theorem~\ref{thm:equivalence} once, with~$O(mn)$ arithmetic operations on integers of polynomial bit length (Lemma~\ref{lem:encoding}), run Wolfe's method with the given insertion rule from the corral~$\hat{B}_0$ alongside, and pivot on every index that Wolfe's method inserts.
By Theorem~\ref{thm:equivalence}(c), this is a simplex run from~$B_0$ that has the same length as the Wolfe run, which proves the statement.
\qed
\end{proof}

\subsection{Proof of Corollary~\ref{cor:intro-rules}}
\label{app:cor-rules}

Corollary~\ref{cor:intro-rules} combines two results.
First, every insertion rule in Table~\ref{tab:wolfe-rules} inherits the lower bound from its corresponding simplex pivot rule, when Wolfe's method is started from a suitable corral (Lemma~\ref{lem:rules}).
Second, the exponential lower bounds on the linopt and minnorm rules still hold when Wolfe's method is started from the atom of minimum norm (Lemma~\ref{lem:km-start}).
Both results are not immediate from Theorem~\ref{thm:equivalence}: the lower bounds for the memory-based and randomized simplex rules in Table~\ref{tab:wolfe-rules} are proven on linear programs whose feasible region can be unbounded, so we add an artificial cap constraint to guarantee that assumption~\ref{ass:bounded} is satisfied; the corral that only contains the atom of minimum norm does not correspond to a feasible basis, so we force Wolfe's method to first include the atoms of some suitable corral, which does correspond to a feasible basis, and use the equivalence to the simplex method from there on.

We start with the proof of the first result.

\begin{lemma}\label{lem:rules}
With every insertion rule of Table~\ref{tab:wolfe-rules}, Wolfe's method requires at least as many iterations (in expectation, for randomized rules) in the worst case as listed in Table~\ref{tab:wolfe-rules}, when started from a suitable corral.
\end{lemma}
\begin{proof}
By Theorem~\ref{thm:equivalence}(a), at every corral $\hat{B}$ that corresponds to a feasible basis $B$, the index of the improving atom of minimum gap at $\vec{w_{\hat B}}$ is exactly the index of minimum reduced cost at $B$; if there is a tie, then the set of indices of minimum gap is the set that contains exactly the indices of minimum reduced cost.
Thus, when using the same tie-breaking rule, Dantzig's rule and the linopt rule always choose the same index.
The same holds for the other rules in Table~\ref{tab:wolfe-rules}, since their choices depend only on the set of improving indices and the history of the run, which agree by Theorem~\ref{thm:equivalence}(a)--(c).
For the randomized rules, this holds for every outcome of the random choices, so the runs agree in distribution.

For the lower bounds, it now suffices to give, for every pivot rule of Table~\ref{tab:wolfe-rules}, linear programs that satisfy \ref{ass:integral}--\ref{ass:nondeg} and feasible bases~$B_0$ from which the pivot rule, with a tie-breaking that depends only on the indices of the current basis and the history, attains the listed lower bound.
For Dantzig's rule and Bland's rule (after renumbering the indices), these are the linear programs~\eqref{eq:km} of Klee and Minty with~$2n$ variables, on which both rules perform~$2^n-1$ pivots from~$B_0=\{s_1,\dots,s_n\}$ by Lemma~\ref{lem:km} below.

The lower bounds for the pivot rules of Zadeh~\cite{DFH23} and Cunningham~\cite{AF17}, for Random-Edge~\cite{FHZ11,Han12}, and for Random-Bland~\cite{FHZ14} are proved on linear programs that are induced by Markov decision processes \cite{Put94,Han12}.
These LPs are nondegenerate and their constraint matrices have full row rank~\cite[Section~2.4.1]{Han12}.
Further, we can scale their rational data to integers without changing the simplex runs (see the discussion in Section~\ref{sec:equivalence}).
However, their feasible sets can be unbounded, as in~\cite[Section~5]{DFH23}, so these programs need not satisfy~\ref{ass:bounded}.
Since all of these programs have a finite optimal solution~\cite[proof of Lemma~1]{DFH23}, \cite[Theorem~7]{AF17}, \cite[Lemma~4.6.1 and Theorems~2.4.6 and~2.4.12]{Han12}, \cite[Lemma~5.5]{FHZ14},
we can use the following standard trick, see e.g.~\cite{BT97}, to bound the feasible set without changing the simplex runs.
If~\eqref{eq:lp} has integer data whose absolute values are at most~$U$, then~$\one^\top\ivec{x}{B}\leq m(mU)^m$ for every basic solution~$\ivec{x}{B}$, and we add a variable~$x_{\mathrm{cap}}$ with objective coefficient~$0$ and the constraint~$\one^\top\vec{x}+x_{\mathrm{cap}}=m(mU)^m+1$.
Every feasible basis~$B$ of~\eqref{eq:lp} gives the feasible basis~$B\cup\{\mathrm{cap}\}$ with the same reduced costs, and~$x_{\mathrm{cap}}$ stays positive along every pivot from~$B$, whose step is finite since the optimum is finite; so the simplex runs from~$B_0$ and from~$B_0\cup\{\mathrm{cap}\}$ coincide.
The new feasible set is bounded, and every new feasible basis is nondegenerate, since a degenerate one would give a basic feasible solution of~\eqref{eq:lp} with fewer than~$m$ positive coordinates or with~$\one^\top\vec{x}>m(mU)^m$.
Hence, the capped LPs satisfy \ref{ass:integral}--\ref{ass:nondeg} and have the same simplex runs from their initial bases.
The tie-breaking rule of~\cite{DFH23} depends only on the current policy and the numbers of previous entries, while the other rules do not produce ties, so all lower bounds of Table~\ref{tab:wolfe-rules} apply to the corresponding insertion rules; the cap adds only one variable, which does not change these bounds.
\qed
\end{proof}

We now derive the second result.
Theorem~\ref{thm:equivalence} starts Wolfe's method at a suitable corral of~$m+1$ atoms, whereas De~Loera et al. \cite{DHR} start it at an atom of minimum norm.
We now modify the linear program, without changing its simplex runs, such that Wolfe's method, started at the atom of minimum norm, first inserts the atoms of a suitable corral one by one and without minor cycles, and from there on simulates the simplex method as in Theorem~\ref{thm:equivalence}.

For this, let~\eqref{eq:lp} denote a linear program satisfying \ref{ass:integral}--\ref{ass:nondeg}, let~$B_0$ be a feasible basis, and let~$\Delta_0=|\det A_{B_0}|$.
We first bring~\eqref{eq:lp} into canonical form with respect to~$B_0$, as in the full tableau implementation of the simplex method~\cite[Section~3.3]{BT97}, but scaled to integer data:
we replace the integer data~$A$,~$\vec{b}$,~$\vec{c}$ of~\eqref{eq:lp} by
\[
    A'\coloneqq\Delta_0A_{B_0}^{-1}A,\qquad\vec{b'}\coloneqq2\Delta_0A_{B_0}^{-1}\vec{b},\qquad\vec{c'}\coloneqq\Delta_0\bigl(\vec{c}-A^\top\ivec{y}{B_0}\bigr).
\]
These new data are integral as well, since~$\Delta_0A_{B_0}^{-1}$ is the adjugate of~$A_{B_0}$ up to sign.
Moreover, we have~$A'_{B_0}=\Delta_0I$ and~$\vec{c'_{B_0}}=\zero$, and the vector $\vec{b'}=2\Delta_0\ivec{x}{B_0}_{B_0}$ has positive even entries.
For convenience, we index the rows of~$A'$ by the indices in~$B_0$, so that we have~$A'_{ii}=\Delta_0$ for~$i\in B_0$.
We call the linear program $\min\{(\vec{c'})^\top\vec{x}:A'\vec{x}=\vec{b'},\ \vec{x}\geq\zero\}$ the \emph{canonical LP}.
Then, we add a new variable~$x_\omega$, the \emph{shield}, with objective coefficient~$0$, which does not occur in the~$m$ constraints of the canonical LP, and add the new constraint
\begin{equation}
    \Delta_0x_\omega-a\sum_{k\in[n]\setminus B_0}x_k=1,\qquad\text{where}\quad
    a\coloneqq\Delta_0+\max_{k\in[n]\setminus B_0}\ \sum_{i\in B_0}b'_i\max\{A'_{ik},0\},
    \label{eq:shield}
\end{equation}
whose row in the constraint matrix we index by~$\omega$.\footnote{The set~$[n]\setminus B_0$ can be empty when~$m=n$. Here and in the definition of~$W$ above equation~\eqref{eq:weights}, a maximum over an empty set is taken to be zero.}
We denote the resulting LP by (LP$^+$) and its data by~$A^+$,~$\vec{b^+}$,~$\vec{c^+}$.
That is, the matrix~$A^+$ arises from~$A'$ by appending the column~$\zero$ of the new index~$\omega=n+1$ and then the row of~\eqref{eq:shield}, and $\vec{b^+}=((\vec{b'})^\top,1)^\top$ and $\vec{c^+}=((\vec{c'})^\top,0)^\top$.
So~$A^+$ has~$m+1$ rows indexed by~$B_0\cup\{\omega\}$, and (LP$^+$) has the indices~$[n+1]$. 
In particular,~$A^+_{B_0\cup\{\omega\}}=\Delta_0I$.

\begin{lemma}\label{lem:shield-lp}
The linear program (LP$^+$) defined above satisfies \ref{ass:integral}--\ref{ass:nondeg}, and its data have polynomial bit length.
Its feasible bases are the sets~$B\cup\{\omega\}$ for the feasible bases~$B$ of~\eqref{eq:lp}.
The reduced costs of the indices in~$[n]$ at~$B\cup\{\omega\}$ are~$\Delta_0$ times those at~$B$, and for every improving index, the pivots at~$B$ and at~$B\cup\{\omega\}$ have the same leaving index.
\end{lemma}
\begin{proof}
We first compare the canonical LP with~\eqref{eq:lp}.
Let~$G\coloneqq\Delta_0A_{B_0}^{-1}$, so~$A'=GA$ and~$\vec{b'}=2G\vec{b}$.
Since~$G$ is nonsingular and~$A'_B=GA_B$ for every~$B\subseteq[n]$, a set~$B$ is a basis of the canonical LP if and only if it is a basis of~\eqref{eq:lp}, and~$A'\vec{x}=\vec{b'}$ holds if and only if~$A\vec{x}=2\vec{b}$.
So the feasible set of the canonical LP is~$\mathcal{P}'\coloneqq\{2\vec{x}:\vec{x}\in\polytope\}$.
For every basis~$B$ and every index~$j\in[n]$, we have
\begin{equation}
    (A'_B)^{-1}\vec{b'}=A_B^{-1}G^{-1}\,2G\vec{b}=2A_B^{-1}\vec{b},
    \quad
    (A'_B)^{-1}\vec{A'_j}=A_B^{-1}G^{-1}G\vec{A_j}=A_B^{-1}\vec{A_j}.
    \label{eq:canonical-basis}
\end{equation}
The basis~$B$ therefore has the basic solution $2\ivec{x}{B}$ in the canonical LP, where $\ivec{x}{B}$ denotes the basic solution of~$B$ in~\eqref{eq:lp}.
Further, the basis~$B$ has the same basic directions $\ivec{d}{B,j}$ in the canonical LP as in~\eqref{eq:lp} by Definition~\ref{def:basis_d}.
So feasibility and nondegeneracy of the bases, and boundedness of the LP are preserved, all ratios in~\eqref{eq:ratio} are doubled, and the ratio-test winners are the same.
In particular, by~\ref{ass:bounded}, the feasible set~$\mathcal{P}'$ of the canonical LP is nonempty and bounded.
Let $\vec{y}\coloneqq\ivec{y}{B_0}$.
By~\eqref{eq:canonical-basis} and the definition of~$\vec{c'}$, the reduced cost of~$j\in[n]$ at~$B$ in the canonical LP is
\begin{equation}
\begin{split}
    c'_j-(\vec{c'_B})^\top(A'_B)^{-1}\vec{A'_j}
    &=\Delta_0\bigl(c_j-\vec{y}^\top\vec{A_j}\bigr)-\Delta_0\bigl(\vec{c_B}-A_B^\top\vec{y}\bigr)^\top A_B^{-1}\vec{A_j}\\
    &=\Delta_0\bigl(c_j-\vec{c_B}^\top A_B^{-1}\vec{A_j}\bigr)
    =\Delta_0\rc^B_j.
\end{split}
\label{eq:canonical-rc}
\end{equation}
The matrix~$G$ is the adjugate of~$A_{B_0}$ up to sign, and $\Delta_0\vec{y}=G^\top\vec{c_{B_0}}$, so the data~$A'$,~$\vec{b'}$,~$\vec{c'}$ of the canonical LP are integral, and by Hadamard's inequality, they and the number~$a$ in~\eqref{eq:shield} have polynomial bit length.
Thus the data of (LP$^+$) have polynomial bit length, since their entries are entries of~$A'$,~$\vec{b'}$, and~$\vec{c'}$, or one of the numbers~$0$,~$1$,~$\Delta_0$, and~$-a$.
Moreover, $\vec{c'_{B_0}}=\Delta_0(\vec{c_{B_0}}-A_{B_0}^\top\vec{y})=\zero$ and $\vec{b'}=2\Delta_0\ivec{x}{B_0}_{B_0}>\zero$ by \ref{ass:nondeg}.

We now transfer our insights from the canonical LP to (LP$^+$), starting with its feasible set.
Let~$\varphi\colon\R^n\to\R^{n+1}$ be the affine map
\begin{equation}
    \varphi(\vec{x})\coloneqq\Bigl(\vec{x}^\top,\ \frac{1}{\Delta_0}\Bigl(1+a\sum\nolimits_{k\in[n]\setminus B_0}x_k\Bigr)\Bigr)^\top.
    \label{eq:shield-map}
\end{equation}
The variable~$x_\omega$ occurs only in~\eqref{eq:shield}, so a vector $(\vec{x}^\top,x_\omega)^\top$ satisfies $A^+(\vec{x}^\top,x_\omega)^\top=\vec{b^+}$ if and only if~$A'\vec{x}=\vec{b'}$ and $(\vec{x}^\top,x_\omega)^\top=\varphi(\vec{x})$.
For~$\vec{x}\geq\zero$, the last coordinate of~$\varphi(\vec{x})$ is at least~$1/\Delta_0$, since~$a>0$ by~$\vec{b'}>\zero$.
\footnote{Apart from its bit length, this proof uses no other property of~$a$, so the lemma, except for the bit length, and its proof remain valid if, for every~$k\in[n]\setminus B_0$, the number~$a$ is replaced by a positive integer~$a_k$ in the coefficients of~$x_k$ in~\eqref{eq:shield} and~\eqref{eq:shield-map}; we use this in the proof of Lemma~\ref{lem:weighted}.}
Hence, the feasible set of (LP$^+$) is~$\varphi(\mathcal{P}')$, which is nonempty and bounded since both hold for~$\mathcal{P}'$ and~$\varphi$ is affine, and
\begin{equation}
    x_\omega\geq\frac{1}{\Delta_0}>0\qquad\text{for every feasible point~$(\vec{x}^\top,x_\omega)^\top$ of (LP$^+$)}.
    \label{eq:shield-positive}
\end{equation}

Next, we determine the feasible bases of (LP$^+$) and verify~\ref{ass:integral}--\ref{ass:nondeg}.
By expanding the determinant along the column of~$\omega$, whose only nonzero entry is~$\Delta_0$, we get~$\det A^+_{B\cup\{\omega\}}=\pm\Delta_0\det A'_B$ for every set~$B\subseteq[n]$ of size~$m$.
Hence~$B\cup\{\omega\}$ is a basis of (LP$^+$) if and only if~$B$ is a basis of~\eqref{eq:lp}.
For such~$B$, the first~$m$ rows of~$A^+_{B\cup\{\omega\}}\vec{z}=\vec{b^+}$ read~$A'_B\vec{z_B}=\vec{b'}$, so by~\eqref{eq:canonical-basis} and the characterization of~$A^+\vec{x}=\vec{b^+}$ below \eqref{eq:shield-map}, the basic solution of~$B\cup\{\omega\}$ is
\begin{equation}
    \vec{x^+}=\varphi\bigl(2\ivec{x}{B}\bigr).
    \label{eq:shield-solution}
\end{equation}
If~$B$ is feasible, then $2\ivec{x}{B}\in\mathcal{P}'$, so~$\vec{x^+}$ is feasible, and it is nondegenerate, since~$2\ivec{x}{B}_B>\zero$ by~\ref{ass:nondeg} and~$x^+_\omega>0$ by~\eqref{eq:shield-positive}.
Conversely, if~$\vec{x^+}\geq\zero$, then $2\ivec{x}{B}\geq\zero$ by~\eqref{eq:shield-solution}, so~$B$ is feasible.
By~\eqref{eq:shield-positive}, the index~$\omega$ is basic in every feasible basis of (LP$^+$).
So the feasible bases of (LP$^+$) are the sets~$B\cup\{\omega\}$ for the feasible bases~$B$ of~\eqref{eq:lp}, and they are nondegenerate.
Finally,~$A^+$ has full row rank, since~$A^+_{B_0\cup\{\omega\}}=\Delta_0I$, and its data are integral, since those of the canonical LP are.
Thus (LP$^+$) satisfies \ref{ass:integral}--\ref{ass:nondeg}.

It remains to compare the reduced costs and the leaving indices.
Let~$B$ be a basis of~\eqref{eq:lp}, and let $\vec{y^+}=((\vec{y'})^\top,y_\omega)^\top$ be the dual vector of~$B\cup\{\omega\}$ in (LP$^+$), where~$\vec{y'}$ has the entries of the rows~$B_0$.
By Definition~\ref{def:basis_a}, it solves $(A^+_{B\cup\{\omega\}})^\top\vec{y^+}=\vec{c^+_{B\cup\{\omega\}}}$.
The column of~$\omega$ in~$A^+$ is~$\Delta_0\vec{e_\omega}$, and the column of~$i\in[n]$ is $((\vec{A'_i})^\top,A^+_{\omega i})^\top$, so this system reads
\begin{equation}
    \Delta_0y_\omega=c^+_\omega=0\qquad\text{and}\qquad
    (\vec{A'_i})^\top\vec{y'}+A^+_{\omega i}y_\omega=c'_i\quad\text{for all }i\in B.
    \label{eq:shield-dual}
\end{equation}
Hence~$y_\omega=0$ and $(A'_B)^\top\vec{y'}=\vec{c'_B}$, that is,~$\vec{y'}$ is the dual vector of~$B$ in the canonical LP.
For~$j\in[n]$, we have~$c^+_j=c'_j$, so by~\eqref{eq:canonical-rc}, the reduced cost of~$j$ at~$B\cup\{\omega\}$ is
\[
    c^+_j-(\vec{y^+})^\top\vec{A^+_j}=c'_j-(\vec{y'})^\top\vec{A'_j}=c'_j-(\vec{c'_B})^\top(A'_B)^{-1}\vec{A'_j}=\Delta_0\rc^B_j.
\]

For the leaving indices, let~$B$ be a feasible basis of~\eqref{eq:lp}, let~$j$ be a nonbasic improving index, and let~$\vec{d^+}$ be its basic direction at the basis~$B\cup\{\omega\}$ of (LP$^+$).
By Definition~\ref{def:basis_d},~$A^+\vec{d^+}=\zero$, and the first~$m$ rows give~$A'_B\vec{d^+_B}=-\vec{A'_j}$.
Thus, on~$[n]$, the direction~$\vec{d^+}$ agrees with the basic direction $\ivec{d}{B,j}$ of~$B$ by~\eqref{eq:canonical-basis}.
For~$r\geq0$, the point~$\vec{x^+}+r\vec{d^+}$ solves~$A^+\vec{x}=\vec{b^+}$, so it equals $\varphi(2\ivec{x}{B}+r\ivec{d}{B,j})$ by~\eqref{eq:shield-solution} and the characterization above.
Since the last coordinate of~$\varphi(\vec{x})$ is positive for~$\vec{x}\geq\zero$, this point is nonnegative if and only if $2\ivec{x}{B}+r\ivec{d}{B,j}\geq\zero$, and then its coordinate~$\omega$ is positive.
Hence the ratio test at~$B\cup\{\omega\}$ gives the step length~$2\theta^{B,j}$, the index~$\omega$ does not attain its minimum, and the index in~$B$ that attains it is the one that attains the minimum in~\eqref{eq:ratio} at~$B$.
So the pivots on~$j$ at~$B$ and at~$B\cup\{\omega\}$ have the same leaving index.
\qed
\end{proof}

We apply the construction of Section~\ref{sec:equivalence} to the program (LP$^+$), which is admissible by Lemma~\ref{lem:shield-lp}.
This gives~$n+2$ integer atoms $\ivec{p}{j}\in\Z^{m+2}$,
$j\in[n+2]$, with the slack index~$n+2$.
Explicitly, let~$T$,~$C$, and~$N$ be the numbers in~\eqref{eq:T}, \eqref{eq:gamma}, and~\eqref{eq:scale} for (LP$^+$), and let~$\vec{\beta}\coloneqq\vec{b^+}/T$.
By~\eqref{eq:alpha}, \eqref{eq:gamma}, and~\eqref{eq:atoms}, for $j\in[n+1]$, the atoms are
\begin{equation}
    \ivec{p}{j}=\bigl(N(C+c^+_j),\,N^2(\vec{A^+_j}-\vec{\beta})^\top\bigr)^\top,\qquad
    \ivec{p}{n+2}=\bigl(NC,\,-N^2\vec{\beta}^\top\bigr)^\top,
    \label{eq:shield-atoms}
\end{equation}
where the constraint coordinates are indexed by the rows~$B_0\cup\{\omega\}$ of~$A^+$.
The corral of the basis~$B_0\cup\{\omega\}$ is~$B_0\cup\{\omega,n+2\}$.

\begin{lemma}\label{lem:init}
Consider the atoms~\eqref{eq:shield-atoms}, and a set~$S$ with~$n+2\in S\subseteq B_0\cup\{\omega,n+2\}$.
Then,
\begin{enumerate}
    \item[(a)] the set~$S$ is a corral,
    \item[(b)] if~$\omega\notin S$, then $\Gamma(\vec{w_S})=(B_0\cup\{\omega\})\setminus S$, with the gaps~$g_\ell(\vec{w_S})=-N^4\Delta_0\beta_\ell$ for~$\ell\in\Gamma(\vec{w_S})$, and every major cycle of Wolfe's method at the corral~$S$ that inserts an index~$\ell\in\Gamma(\vec{w_S})$ consists of no minor cycle and ends at the corral~$S\cup\{\ell\}$,
    \item[(c)] the atom $\ivec{p}{n+2}$ is the unique atom of minimum norm.
\end{enumerate}
\end{lemma}

\begin{proof}
By~\eqref{eq:shield-atoms}, since~$A^+_{B_0\cup\{\omega\}}=\Delta_0I$ and~$\vec{c^+_{B_0\cup\{\omega\}}}=\zero$, the atoms of the indices in~$B_0\cup\{\omega,n+2\}$ are
\begin{equation}
    \ivec{p}{\ell}=\bigl(NC,\,N^2(\Delta_0\vec{e_\ell}-\vec{\beta})^\top\bigr)^\top,\qquad
    \ivec{p}{n+2}=\bigl(NC,\,-N^2\vec{\beta}^\top\bigr)^\top,
    \label{eq:init-atoms}
\end{equation}
for $\ell\in B_0\cup\{\omega\}$.
The vector~$\vec{\beta}=\vec{b^+}/T$ is positive, as~$\vec{b^+}=((\vec{b'})^\top,1)^\top$ and~$\vec{b'}>\zero$ by construction.
The basic solution~$\vec{x^+}$ of the basis~$B_0\cup\{\omega\}$ of the program (LP$^+$) has the entries~$x^+_\ell=b^+_\ell/\Delta_0$ for~$\ell\in B_0\cup\{\omega\}$, so Lemma~\ref{lem:T-bound} gives
\begin{equation}
    \sum_{\ell\in B_0\cup\{\omega\}}\beta_\ell=\frac{\Delta_0}{T}\,\one^\top\vec{x^+}<\frac{\Delta_0}{2}.
    \label{eq:beta-sum}
\end{equation}

\emph{Part (a).} Let~$R\coloneqq(B_0\cup\{\omega\})\setminus S$, and let~$\vec{\beta_R}$ denote the vector~$\vec{\beta}$ with the entries outside~$R$ replaced by zero.
For affine coefficients~$\vec{\lambda}$ on~$S$, the atoms~\eqref{eq:init-atoms} give~$\sum_{i\in S}\lambda_i\ivec{p}{i}=\bigl(NC,\,N^2(\Delta_0\vec{\lambda_{S\setminus\{n+2\}}}-\vec{\beta})^\top\bigr)^\top$, where~$\vec{\lambda_{S\setminus\{n+2\}}}$ is extended by zero to~$B_0\cup\{\omega\}$.
This sum determines~$\vec{\lambda}$, so the atoms of~$S$ are affinely independent, and its norm is minimized exactly by~$\lambda_\ell=\beta_\ell/\Delta_0$ for~$\ell\in S\setminus\{n+2\}$, which gives
\[
    \vec{w_S}=\bigl(NC,\,-N^2\vec{\beta_R}^\top\bigr)^\top,\qquad
    \icoord{\bar{\lambda}}{S}{{n+2}}=1-\sum_{\ell\in S\setminus\{n+2\}}\frac{\beta_\ell}{\Delta_0}>\frac12
\]
by~\eqref{eq:beta-sum} and~$\vec{\beta}>\zero$.
All coefficients are positive, so~$S$ is a corral.

\emph{Part (b).} Let~$R$ and~$\vec{\beta_R}$ be as in part~(a).
We have~$\norm{\vec{w_S}}^2=N^2C^2+N^4\norm{\vec{\beta_R}}^2$, and by~\eqref{eq:init-atoms}, every~$\ell\in B_0\cup\{\omega\}$ satisfies
\[
    \vec{w_S}^\top\ivec{p}{\ell}=
    \begin{cases}
        N^2C^2+N^4\norm{\vec{\beta_R}}^2-N^4\Delta_0\beta_\ell, & \text{if }\ell\in R,\\
        N^2C^2+N^4\norm{\vec{\beta_R}}^2, & \text{if }\ell\notin R,
    \end{cases}
\]
so~$g_\ell(\vec{w_S})=-N^4\Delta_0\beta_\ell<0$ for~$\ell\in R$, and the gaps of the indices in~$S$ vanish.
For~$k\in[n]\setminus B_0$, the same computation with the atom $\ivec{p}{k}$ in~\eqref{eq:shield-atoms} gives
\[
    g_k(\vec{w_S})=N^2Cc^+_k-N^4\sum_{\ell\in R}A^+_{\ell k}\beta_\ell.
\]
Since~$\omega\notin S$, we have~$\omega\in R$.
By~\eqref{eq:shield}, we have~$A^+_{\omega k}=-a$ and~$b^+_\omega=1$, so the choice of~$a$ there gives
\[
    \sum_{\ell\in R}A^+_{\ell k}\beta_\ell\leq\frac1T\Bigl(-a+\sum_{i\in B_0}b^+_i\max\{A^+_{ik},0\}\Bigr)<0.
\]
Since~$A^+$ and~$\vec{b^+}=T\vec{\beta}$ are integral, the left-hand side of this inequality lies in~$\frac1T\Z$, and hence it is at most~$-\frac1T$,
and~$g_k(\vec{w_S})\geq N^4/T-N^2C|c^+_k|>0$, since~$N^2>TK^2$ by~\eqref{eq:scale} and~$|c^+_k|<C\leq K$ by~\eqref{eq:gamma}.
So $\Gamma(\vec{w_S})=R$.
If a major cycle of Wolfe's method at~$S$ inserts~$\ell\in R$, then the new support~$S\cup\{\ell\}$ is a corral by~(a), so the loop condition of the minor cycles in Line~\ref{line:wf-loop} of Algorithm~\ref{alg:wolfe} fails, and the major cycle ends at the corral~$S\cup\{\ell\}$ without a minor cycle.

\emph{Part (c).} For~$\ell\in B_0\cup\{\omega\}$, we have $\norm{\ivec{p}{\ell}}^2-\norm{\ivec{p}{n+2}}^2=N^4\Delta_0(\Delta_0-2\beta_\ell)$ by~\eqref{eq:init-atoms}, which is positive by~\eqref{eq:beta-sum} and~$\vec{\beta}>\zero$.
For~$k\in[n]\setminus B_0$,
\[
    \norm{\ivec{p}{k}}^2-\norm{\ivec{p}{n+2}}^2
    =N^2\bigl((C+c^+_k)^2-C^2\bigr)+N^4\sum_{\ell\in B_0\cup\{\omega\}}A^+_{\ell k}\bigl(A^+_{\ell k}-2\beta_\ell\bigr).
\]
Every term of the sum with~$\ell\in B_0$ is at least~$-\beta_\ell^2$, since it equals~$(A^+_{\ell k}-\beta_\ell)^2-\beta_\ell^2$,
and the term of~$\omega$ is~$a(a+2\beta_\omega)\geq a^2\geq\Delta_0^2$.
By~\eqref{eq:beta-sum} and~$\vec{\beta}>\zero$, we have~$\sum_{\ell\in B_0}\beta_\ell^2<\Delta_0^2/4$, so the sum is greater than~$\frac34\Delta_0^2\geq\frac34$.
The first term is at least~$-N^2C^2$, and $N>2C$ by~\eqref{eq:scale}, since~$H\geq K\geq C$, so the difference is positive, and~$\ivec{p}{n+2}$ is the unique atom of minimum norm.
\qed
\end{proof}

This allows us to relate the runs of Wolfe's method with the linopt and least-index rules to the simplex runs.
The key argument is that both insertion rules choose $\omega$ only after the indices in $B_0$.
The same does not necessarily hold for the minnorm rule, which is the least-index rule for the order of the atoms by norm.
\begin{lemma}\label{lem:init-run}
The atoms~\eqref{eq:shield-atoms} can be computed in polynomial time.
From the unique shortest atom $\ivec{p}{n+2}$, Wolfe's method with the linopt rule, or the least-index rule for an order with~$\omega$ after~$[n]$, first inserts the atoms of~$B_0\cup\{\omega\}$ one-by-one without minor cycles; then its corrals are the sets~$B_k\cup\{\omega,n+2\}$,~$k=0,\dots,f$, for a simplex run~$B_0,\dots,B_f$ on~\eqref{eq:lp} with Dantzig's or Bland's rule, respectively.
\end{lemma}

\begin{proof}
By Lemma~\ref{lem:shield-lp}, (LP$^+$) satisfies \ref{ass:integral}--\ref{ass:nondeg}, so Section~\ref{sec:equivalence} applies to it.
Its data have polynomial bit length by Lemma~\ref{lem:shield-lp}, so Gaussian elimination computes them in polynomial time, and Lemma~\ref{lem:encoding} gives the atoms in polynomial time.
By Lemma~\ref{lem:init}(c), the atom $\ivec{p}{n+2}$ is the unique atom of minimum norm, so the method starts at the corral~$\{n+2\}$.
As long as the current corral~$S\subseteq B_0\cup\{\omega,n+2\}$ does not contain~$\omega$, Lemma~\ref{lem:init}(b) shows that $\Gamma(\vec{w_S})=(B_0\cup\{\omega\})\setminus S$ and that the next corral is~$S\cup\{\ell\}$ for the inserted index~$\ell$, without a minor cycle.
By Lemma~\ref{lem:init}(b), the gaps of these indices are~$g_\ell(\vec{w_S})=-N^4\Delta_0\beta_\ell$, so the linopt rule, which minimizes the gap, inserts an index~$\ell$ that maximizes~$\beta_\ell=b^+_\ell/T$.
Recall that we have~$\vec{b^+}=((\vec{b'})^\top,1)^\top$, where the vector~$\vec{b'}=2\Delta_0\ivec{x}{B_0}_{B_0}$ has positive even entries.
Hence~$b^+_i=b'_i\geq2>1=b^+_\omega$ for~$i\in B_0$, and the linopt rule inserts~$\omega$ only after all indices of~$B_0$.
The least-index rule does the same by the choice of the order.
So after~$m+1$ major cycles, the method is at the corral~$B_0\cup\{\omega,n+2\}$ of the basis~$B_0\cup\{\omega\}$, with the coefficients of its affine minimizer.
From there, Theorem~\ref{thm:equivalence}(c) and the proof of Lemma~\ref{lem:rules} for (LP$^+$) show that the method visits the corrals~$B^+_k\cup\{n+2\}$,~$k=0,\dots,f$, for a simplex run~$B^+_0,\dots,B^+_f$ on (LP$^+$) from~$B^+_0=B_0\cup\{\omega\}$ with Dantzig's or Bland's rule, respectively.
By Lemma~\ref{lem:shield-lp}, the feasible bases of (LP$^+$) are the sets~$B\cup\{\omega\}$ for the feasible bases~$B$ of~\eqref{eq:lp}, the reduced costs of the indices in~$[n]$ at~$B\cup\{\omega\}$ are~$\Delta_0>0$ times those at~$B$, and the pivots at~$B$ and at~$B\cup\{\omega\}$ have the same leaving indices.
Moreover,~$\omega$ is basic in every feasible basis of (LP$^+$), hence never improving.
So~$B^+_k=B_k\cup\{\omega\}$ for a simplex run~$B_0,\dots,B_f$ on~\eqref{eq:lp} with the same rule.
\qed
\end{proof}

It is well-known that the Klee--Minty cube yields an exponential lower bound for Dantzig's~\cite{KM72} and Bland's~\cite{AC78} pivot rules.
For~$n\geq1$, the linear program of Klee and Minty~\cite{KM72}, in the integer form of Chv\'atal~\cite[Chapter~4]{Chv83}, has the~$2n$ variables~$x_1,\dots,x_n$ and~$s_1,\dots,s_n$ and is given by
\begin{equation}
    \min\ -\sum_{j=1}^n10^{n-j}x_j\quad\text{s.t.}\quad
    x_i+s_i+2\sum_{j<i}10^{i-j}x_j=100^{i-1}\ \ (i\in[n]),\quad\vec{x},\vec{s}\geq\zero.
    \label{eq:km}
\end{equation}

It is folklore that its feasible set is nonempty and bounded and that its feasible bases are nondegenerate (see, e.g.,~\cite[Section~3]{AC78}), so~\eqref{eq:km} satisfies \ref{ass:integral}--\ref{ass:nondeg}.
We are not aware of a lower bound for Bland's rule on this specific formulation that uses the order in Lemma~\ref{lem:km}, and that states the uniqueness of the choices of Dantzig's rule, so we include a short proof for both.

\begin{lemma}\label{lem:km}
On the linear program~\eqref{eq:km}, started at the basis~$B_0\coloneqq\{s_1,\dots,s_n\}$, Dantzig's rule and Bland's rule for the order~$x_1\prec s_1\prec x_2\prec s_2\prec\dots\prec x_n\prec s_n$ perform~$2^n-1$ pivots, and every choice of Dantzig's rule is unique.
\end{lemma}

\begin{proof}
Substituting~$x_j=100^{j-1}\tilde{x}_j$ and~$s_i=100^{i-1}\tilde{s}_i$ and dividing row~$i$ by~$100^{i-1}$ and the objective by~$10^{2n-2}$ turns~\eqref{eq:km} into the example~(3) of Avis and Chv\'atal~\cite[Section~3]{AC78} with~$\varepsilon=1/10$.
This rescaling of the variables by positive factors changes neither the feasible bases, nor the signs of the reduced costs, nor the leaving indices (see the proof of Lemma~\ref{lem:weighted}(a)).
By~\cite[Section~3]{AC78}, the feasible bases are the sets that contain exactly one index of every pair~$(x_i,s_i)$.
We encode such a basis by the vector~$\vec{a}\in\{0,1\}^n$ with~$a_i=1$ if~$x_i$ is basic and~$a_i=0$ if~$s_i$ is basic.
Moreover, by~\cite[Section~3]{AC78}, at the feasible basis with the vector~$\vec{a}$, the nonbasic index of the pair~$(x_k,s_k)$ is improving if and only if~$a_k+\dots+a_n$ is even, and the pivot on the improving pair with the smallest indices leads to the immediate successor of~$\vec{a}$ in a linear order of~$\{0,1\}^n$ whose first element is~$(0,\dots,0)$, the vector of~$B_0$, and whose last element is~$(0,\dots,0,1)$, the vector of the optimal basis.
Bland's rule for the given order chooses the improving pair with the smallest indices, so it visits all~$2^n$ feasible bases and performs~$2^n-1$ pivots.
Dantzig's rule makes the same choices, uniquely, since the reduced cost of the nonbasic variable~$w_j$ of the pair whose indices are~$j$ has the absolute value~$10^{n-j}$, which decreases with~$j$.
To see this, let~$Z_i\coloneqq\sum_{j\leq i}10^{i-j}x_j$ for~$i\in[n]$ and~$Z_0\coloneqq0$.
Row~$i$ reads~$x_i+s_i+20Z_{i-1}=100^{i-1}$, and~$Z_i=x_i+10Z_{i-1}$, so~$Z_i=w_i+10Z_{i-1}$ if~$a_i=0$, and~$Z_i=100^{i-1}-w_i-10Z_{i-1}$ if~$a_i=1$.
By induction on~$i$, there are numbers~$\kappa_0,\dots,\kappa_n$ that depend only on~$\vec a$ such that every solution of the constraints of~\eqref{eq:km} satisfies
\[
    Z_i=\kappa_i+\sum_{j\in[i]}(-1)^{a_j+\dots+a_i}\,10^{i-j}w_j\qquad\text{for all }i\in\{0,\dots,n\}.
\]
For~$i=n$, this writes the objective~$-Z_n$ as an affine function of the nonbasic variables, and by Lemma~\ref{lem:identities}(c), the coefficient~$-(-1)^{a_j+\dots+a_n}10^{n-j}$ of~$w_j$ is the reduced cost of~$w_j$.
Thus, Dantzig's rule also performs~$2^n-1$ pivots, as was already shown in~\cite{KM72} and~\cite[Chapter~4]{Chv83}.
\qed
\end{proof}

Together with Lemma~\ref{lem:init-run}, this lower bound for Dantzig's rule already yields the lower bound for Wolfe's linopt rule from the unique atom of minimum norm (Lemma~\ref{lem:km-start}).
It does not yet give the lower bound for minnorm, since, on the atoms~\eqref{eq:shield-atoms} for~\eqref{eq:km}, the order of the atoms by norm does not agree with the order of Lemma~\ref{lem:km}.
To overcome this, we scale the columns of (LP$^+$) such that the norms of the atoms increase along a given order, which does not change the runs of Bland's rule.

For positive integers~$d_\ell$, with~$\ell\in B_0\cup\{\omega\}$, and integers~$a_k\geq a$, with~$k\in[n]\setminus B_0$, let the linear program ($\widetilde{\text{LP}}{}^+$) arise from (LP$^+$) by replacing the column~$\Delta_0\vec{e_\ell}$ of each index~$\ell\in B_0\cup\{\omega\}$ by~$d_\ell\vec{e_\ell}$, and the coefficient~$(-a)$ of each~$x_k$ in~\eqref{eq:shield} by~$(-a_k)$.
For~$\ell\in B_0$, this substitutes~$(d_\ell/\Delta_0)x_\ell$ for~$x_\ell$, which multiplies the reduced cost of~$\ell$ by~$d_\ell/\Delta_0$ and changes neither the improving sets nor the leaving indices.
In particular, it may change the choices of Dantzig's rule, but not those of Bland's rule.
Let~$\prec$ be an order of~$[n]$, let~$r(\ell)$ be the position of~$\ell\in[n]$ in~$\prec$, let~$r(\omega)\coloneqq n+1$, and let~$W\coloneqq a+\max_{k\in[n]\setminus B_0}\norm{\vec{A'_k}}^2$.
We choose the weights
\begin{equation}
    d_\ell\coloneqq r(\ell)W\quad\text{and}\quad a_k\coloneqq r(k)W,
    \label{eq:weights}
\end{equation}
for $\ell\in B_0\cup\{\omega\}$ and $k\in[n]\setminus B_0$.
Let~$\tilde{A}$ denote the constraint matrix of ($\widetilde{\text{LP}}{}^+$), whose right-hand side is~$\vec{b^+}$ and objective vector is~$\vec{c^+}$.
For~$\ell\in B_0\cup\{\omega\}$, the column~$\vec{\tilde{A}_\ell}=d_\ell\vec{e_\ell}$ has the squared norm~$r(\ell)^2W^2$.
For~$k\in[n]\setminus B_0$, the column~$\vec{\tilde{A}_k}=((\vec{A'_k})^\top,-a_k)^\top$ has the squared norm~$r(k)^2W^2+\norm{\vec{A'_k}}^2$, where~$\norm{\vec{A'_k}}^2\leq W-a<W$ by the choice of~$W$.
So the squared norm of the column of every index~$\ell\in[n]\cup\{\omega\}$ lies in~$[r(\ell)^2W^2,\,r(\ell)^2W^2+W)$, and the column norms increase along~$\prec$, followed by~$\omega$.

We apply the construction of Section~\ref{sec:equivalence} to ($\widetilde{\text{LP}}{}^+$) with the weights~\eqref{eq:weights}, which gives~$n+2$ integer atoms~$\ivec{p}{j}\in\Z^{m+2}$, for~$j\in[n+2]$, with the slack index~$n+2$.
Explicitly, let~$T$,~$C$, and~$N$ now be the numbers in~\eqref{eq:T}, \eqref{eq:gamma}, and~\eqref{eq:scale} for ($\widetilde{\text{LP}}{}^+$), and let~$\vec{\beta}\coloneqq\vec{b^+}/T$.
Similar to~\eqref{eq:shield-atoms}, for $j\in[n+1]$, the atoms are
\begin{equation}
    \ivec{p}{j}=\bigl(N(C+c^+_j),\,N^2(\vec{\tilde{A}_j}-\vec{\beta})^\top\bigr)^\top,\qquad
    \ivec{p}{n+2}=\bigl(NC,\,-N^2\vec{\beta}^\top\bigr)^\top.
    \label{eq:weighted-atoms}
\end{equation}

The following lemma repeats the arguments of Lemmas~\ref{lem:shield-lp} and~\ref{lem:init} for these atoms.
This is necessary for our analysis: the previous lemmas do not cover the minnorm rule (see the discussion above), while Lemma~\ref{lem:minnorm} below does not hold for the linopt rule, since the scaling that we use to obtain ($\widetilde{\text{LP}}{}^+$) multiplies the reduced costs of the indices in~$B_0$ by different factors and may therefore change the choices of Dantzig's rule.

\begin{lemma}\label{lem:weighted}
Consider the atoms~\eqref{eq:weighted-atoms}.
\begin{enumerate}
    \item[(a)] The linear program ($\widetilde{\text{LP}}{}^+$) satisfies \ref{ass:integral}--\ref{ass:nondeg}, and its data have polynomial bit length. Its feasible bases are the sets~$B\cup\{\omega\}$ for the feasible bases~$B$ of~\eqref{eq:lp}, its improving indices at~$B\cup\{\omega\}$ are those at~$B$, and for every improving index, the pivots at~$B$ and at~$B\cup\{\omega\}$ have the same leaving index.
    \item[(b)] Let~$S$ be a set with~$n+2\in S\subseteq B_0\cup\{\omega,n+2\}$. Then~$S$ is a corral. If~$\omega\notin S$, then the improving atoms at~$\vec{w_S}$ are exactly those of the indices in~$(B_0\cup\{\omega\})\setminus S$, and every major cycle of Wolfe's method at the corral~$S$ that inserts such an index~$\ell$ consists of no minor cycle and ends at the corral~$S\cup\{\ell\}$.
    \item[(c)] The norms of the atoms increase strictly in the following order of their indices: first the index~$n+2$, then the indices in~$[n]$ ordered by~$\prec$, then the index~$\omega$.
\end{enumerate}
\end{lemma}

\begin{proof}
\emph{Part (a).} A substitution~$x_\ell\mapsto s_\ell x_\ell$ with~$s_\ell>0$ multiplies the column and the objective coefficient of~$\ell$ by~$s_\ell$.
It changes neither the bases nor the dual vectors, since~$(A_BS_B)^{-\top}S_B\vec{c_B}=A_B^{-\top}\vec{c_B}$ for the diagonal matrix~$S_B$ of the factors, and it multiplies the reduced cost of~$\ell$ by~$s_\ell$.
It divides the coordinate~$\ell$ of the basic solution by~$s_\ell$ and multiplies the coordinate~$\ell$ of the basic direction of~$j$ by~$s_j/s_\ell$, so all ratios in~\eqref{eq:ratio} are divided by~$s_j$, and the ratio-test winners stay the same.
Since~$a_k\geq a>0$, the footnote in the proof of Lemma~\ref{lem:shield-lp} shows that Lemma~\ref{lem:shield-lp}, except for the bit length, holds for the linear program that arises from (LP$^+$) by replacing the coefficient~$(-a)$ of each~$x_k$ in~\eqref{eq:shield} by~$(-a_k)$.
The program ($\widetilde{\text{LP}}{}^+$) arises from this linear program by the substitutions~$x_\ell\mapsto(d_\ell/\Delta_0)x_\ell$ for~$\ell\in B_0\cup\{\omega\}$, which replace the columns~$\Delta_0\vec{e_\ell}$ by~$d_\ell\vec{e_\ell}$ and keep the objective coefficients~$c^+_\ell=0$.
Hence both programs have the same feasible bases, the same improving indices at them, and the same leaving indices, and the feasible bases of ($\widetilde{\text{LP}}{}^+$) are nondegenerate.
Moreover, the constraint matrix of ($\widetilde{\text{LP}}{}^+$) has full row rank, its data are integral, and its feasible set is nonempty and bounded, since it arises from the feasible set of the other program by multiplying coordinates by positive factors.
So ($\widetilde{\text{LP}}{}^+$) satisfies \ref{ass:integral}--\ref{ass:nondeg}, its feasible bases are the sets~$B\cup\{\omega\}$ for the feasible bases~$B$ of~\eqref{eq:lp}, its improving indices at~$B\cup\{\omega\}$ are those at~$B$, and its pivots have the same leaving indices.
The data of the canonical LP have polynomial bit length by Lemma~\ref{lem:shield-lp}, and so do~$W$,~$d_\ell$, and~$a_k$; hence so do the data of~($\widetilde{\text{LP}}{}^+$).

\emph{Part (b).} The diagonal matrix~$\tilde{A}_{B_0\cup\{\omega\}}$ has the entries~$d_\ell$, for~$\ell\in B_0\cup\{\omega\}$.
Thus, by~\eqref{eq:weighted-atoms}, we have~$\ivec{p}{\ell}=(NC,\,N^2(d_\ell\vec{e_\ell}-\vec{\beta})^\top)^\top$ for~$\ell\in B_0\cup\{\omega\}$, and the basic solution of the basis~$B_0\cup\{\omega\}$ has the coordinates~$b^+_\ell/d_\ell$.
Lemma~\ref{lem:T-bound} gives~$\sum_{\ell\in B_0\cup\{\omega\}}\beta_\ell/d_\ell<\frac12$, where~$\vec{\beta}=\vec{b^+}/T$ and~$\vec{b^+}=((\vec{b'})^\top,1)^\top$ is the right-hand side of (LP$^+$) and of ($\widetilde{\text{LP}}{}^+$), and the proof of Lemma~\ref{lem:init}(a),(b) applies with~$\tilde{A}$ in place of~$A^+$ and~$\Delta_0$ replaced by~$d_\ell$: the affine minimizer of~$S$ has the coefficients~$\beta_\ell/d_\ell$ for~$\ell\in S\setminus\{n+2\}$, the gaps of~$\ell\in(B_0\cup\{\omega\})\setminus S$ are~$-N^4d_\ell\beta_\ell$, and the estimate for~$k\in[n]\setminus B_0$ uses only~$b^+_\omega=1$,~$\tilde{A}_{\omega k}=-a_k$, and~$a_k\geq a$.
This proves~(b).

\emph{Part (c).} The slack index~$n+2$ has no column in~$\tilde{A}$, and we set~$\vec{\tilde{A}_{n+2}}\coloneqq\zero$ and~$c^+_{n+2}\coloneqq0$, so that~\eqref{eq:weighted-atoms} also gives~$\ivec{p}{n+2}$ by the formula for~$j\in[n+1]$, and let~$\gamma_\ell\coloneqq C+c^+_\ell$ for~$\ell\in[n+2]$.
By~\eqref{eq:weighted-atoms}, every index~$\ell\in[n+2]$ satisfies
\[
    \norm{\ivec{p}{\ell}}^2=N^2\gamma_\ell^2+N^4\bigl(\norm{\vec{\tilde{A}_\ell}}^2-2(\vec{\tilde{A}_\ell})^\top\vec{\beta}+\norm{\vec{\beta}}^2\bigr).
\]
We first show that the terms other than~$N^4\norm{\vec{\tilde{A}_\ell}}^2$ are small.
We have~$1\leq\gamma_\ell<2C$ by~\eqref{eq:gamma}.
Let~$U\geq1$ now be the number in~\eqref{eq:T} for ($\widetilde{\text{LP}}{}^+$), that is, the largest absolute value of an entry of~$\tilde{A}$ or~$\vec{b^+}$.
Since ($\widetilde{\text{LP}}{}^+$) has~$m+1\geq2$ rows, the definition in \eqref{eq:T} gives~$T>2(m+1)\bigl((m+1)U\bigr)^{m+1}\geq2(m+1)^3U^2$, so~$|(\vec{\tilde{A}_\ell})^\top\vec{\beta}|\leq(m+1)U^2/T<\frac18$.
Hence, for all indices~$j,k\in[n+2]$ such that~$\norm{\vec{\tilde{A}_k}}^2\geq\norm{\vec{\tilde{A}_j}}^2+1$, we have
\[
    \norm{\ivec{p}{k}}^2-\norm{\ivec{p}{j}}^2>N^4\Bigl(1-\frac12\Bigr)-4N^2C^2>0,
\]
because~$N>TC\geq3C$ by~\eqref{eq:scale}.
It remains to compare the values~$\norm{\vec{\tilde{A}_\ell}}^2$.
Let~$r(n+2)\coloneqq0$.
We have~$\norm{\vec{\tilde{A}_{n+2}}}^2=0$ by definition, and~$r(\ell)^2W^2\leq\norm{\vec{\tilde{A}_\ell}}^2\leq r(\ell)^2W^2+W-1$ for~$\ell\in[n]\cup\{\omega\}$ by the discussion before the lemma.
Therefore~$r(j)<r(k)$ implies that~$\norm{\vec{\tilde{A}_k}}^2-\norm{\vec{\tilde{A}_j}}^2\geq(r(j)+1)^2W^2-r(j)^2W^2-W+1\geq1$ for all indices~$j,k\in[n+2]$.
This proves~(c).
\qed
\end{proof}

The proof of the following lemma follows the same steps as the proof of Lemma~\ref{lem:init-run}, with Lemma~\ref{lem:weighted} in place of Lemmas~\ref{lem:shield-lp} and~\ref{lem:init}.

\begin{lemma}\label{lem:minnorm}
The atoms~\eqref{eq:weighted-atoms} can be computed in polynomial time and have distinct norms.
From the unique shortest atom $\ivec{p}{n+2}$, Wolfe's method with the minnorm rule first inserts the atoms of~$B_0\cup\{\omega\}$ one-by-one without minor cycles; then its corrals are the sets~$B_k\cup\{\omega,n+2\}$,~$k=0,\dots,f$, for the simplex run~$B_0,\dots,B_f$ on~\eqref{eq:lp} with Bland's rule for~$\prec$.
\end{lemma}

\begin{proof}
By Lemma~\ref{lem:weighted}(a), Section~\ref{sec:equivalence} applies to ($\widetilde{\text{LP}}{}^+$), and its data have polynomial bit length, so they can be computed in polynomial time as in the proof of Lemma~\ref{lem:init-run}.
Thus, Lemma~\ref{lem:encoding} gives the atoms in polynomial time.
By Lemma~\ref{lem:weighted}(c), their norms are pairwise distinct, and $\ivec{p}{n+2}$ is the shortest atom.
So the method starts at the corral~$\{n+2\}$.
As long as the current corral $S\subseteq B_0\cup\{\omega,n+2\}$ does not contain~$\omega$, Lemma~\ref{lem:weighted}(b) shows that the improving atoms are those of the indices in~$(B_0\cup\{\omega\})\setminus S$ and that the next corral is~$S\cup\{\ell\}$ for the inserted index~$\ell$, without a minor cycle.
By Lemma~\ref{lem:weighted}(c), $\ivec{p}{\omega}$ is the longest of these atoms, so the minnorm rule inserts the indices of~$B_0$ before~$\omega$; after~$m+1$ major cycles, the method is at the corral~$B_0\cup\{\omega,n+2\}$ of the basis~$B_0\cup\{\omega\}$, with the coefficients of its affine minimizer.
From there on, the minnorm rule is the least-index rule for the order by norm, whose restriction to~$[n]$ is~$\prec$, since~$\omega$ is basic and the slack index~$n+2$ has gap zero.
Theorem~\ref{thm:equivalence}(c) and the proof of Lemma~\ref{lem:rules} for ($\widetilde{\text{LP}}{}^+$) show that the method visits the corrals~$B^+_k\cup\{n+2\}$,~$k=0,\dots,f$, for the simplex run~$B^+_0,\dots,B^+_f$ on ($\widetilde{\text{LP}}{}^+$) from~$B^+_0=B_0\cup\{\omega\}$ with Bland's rule for~$\prec$.
By Lemma~\ref{lem:weighted}(a),~$B^+_k=B_k\cup\{\omega\}$ for the simplex run~$B_0,\dots,B_f$ on~\eqref{eq:lp} with Bland's rule for~$\prec$.
\qed
\end{proof}

We can now obtain the lower bounds for linopt and minnorm as consequences of the Klee--Minty cube.

\begin{lemma}\label{lem:km-start}
For every~$n\geq1$, there are~$2n+2$ atoms in~$\Z^{n+2}$, with encoding length polynomial in~$n$, and a unique atom of minimum norm, on which Wolfe's method with the linopt rule, started at this atom, visits exactly~$2^n+n+1$ corrals, for every rule for breaking ties; the same holds for the minnorm rule, on other such atoms.
\end{lemma}

\begin{proof}
We first consider the linopt rule.
The Klee--Minty program~\eqref{eq:km} is admissible, i.e., it satisfies \ref{ass:integral}--\ref{ass:nondeg}, as noted below~\eqref{eq:km}.
We take~\eqref{eq:km} as the linear program~\eqref{eq:lp} in the construction of (LP$^+$), so it has~$m=n$ constraints and~$2n$ indices, and we take the basis~$B_0=\{s_1,\dots,s_n\}$ of Lemma~\ref{lem:km}.
Then (LP$^+$) has~$n+1$ constraints and the~$2n+1$ indices~$[2n]\cup\{\omega\}$ with~$\omega=2n+1$, so the atoms~\eqref{eq:shield-atoms} of (LP$^+$) are~$2n+2$ integer atoms in~$\Z^{n+2}$ with the slack index~$2n+2$.
By Lemma~\ref{lem:init-run}, they can be computed in time polynomial in~$n$, and the slack atom~$\ivec{p}{2n+2}$ is the unique atom of minimum norm by Lemma~\ref{lem:init}(c).
By Lemma~\ref{lem:init-run} and Lemma~\ref{lem:km}, the run visits~$1+(n+1)+(2^n-1)=2^n+n+1$ corrals.
The linopt rule never faces a tie.
In the first~$n+1$ major cycles, the gap of an improving index~$\ell$ is~$-N^4\Delta_0b^+_\ell/T$ by Lemma~\ref{lem:init}(b).
Here the columns of~$B_0$ in~\eqref{eq:km} form the identity matrix, so~$\Delta_0=1$ and~$\vec{b'}=2\vec{b}$, and the right-hand sides~$b^+_\ell$ are the numbers~$2\cdot100^{i-1}$,~$i\in[n]$, and~$1$, which are distinct.
Afterwards, Dantzig's choices are unique by Lemma~\ref{lem:km}, and so are those of the linopt rule, since its gaps are positive multiples of the reduced costs of~\eqref{eq:km} by Theorem~\ref{thm:equivalence}(a) and Lemma~\ref{lem:shield-lp}.
On a side note, the same count holds for the least-index rule for the order of the variables in Lemma~\ref{lem:km}, followed by $x_\omega$ and the slack $x_{2n+2}$.

For the minnorm rule, we take the atoms~\eqref{eq:weighted-atoms} of ($\widetilde{\text{LP}}{}^+$) for~\eqref{eq:km}, the basis~$B_0$, and the order~$\prec$ of Lemma~\ref{lem:km}, and argue as for the linopt rule, with Lemma~\ref{lem:minnorm} in place of Lemma~\ref{lem:init-run} and Lemma~\ref{lem:weighted} in place of Lemmas~\ref{lem:shield-lp} and~\ref{lem:init}.
By Lemma~\ref{lem:weighted}(c), the atoms have pairwise distinct norms, and~$\ivec{p}{2n+2}$ is the unique atom of minimum norm, so the minnorm rule faces no ties.
\qed
\end{proof}

We now combine the two main results of this section to prove Corollary~\ref{cor:intro-rules}.

\begin{proof}[of Corollary~\ref{cor:intro-rules}]
For every insertion rule of Table~\ref{tab:wolfe-rules}, the proof of Lemma~\ref{lem:rules} provides linear programs with~$n$ variables and encoding length polynomial in~$n$, on which Wolfe's method, started at a suitable corral of the~$n+1$ atoms of Theorem~\ref{thm:equivalence}, needs at least the number of major cycles listed in Table~\ref{tab:wolfe-rules} (in expectation, for randomized rules).
These atoms have encoding length polynomial in~$n$ by Lemma~\ref{lem:encoding}, and the bounds in Table~\ref{tab:wolfe-rules} are exponential in~$n$ for the linopt, least-index, least-inserted, and round-robin rules, and superpolynomial in~$n$ for the random and the random-order rules.
For the linopt and the minnorm rules, Lemma~\ref{lem:km-start} gives exponential lower bounds when Wolfe's method is started at the atom of minimum norm.
\qed
\end{proof}

\subsection{Proof of Corollary~\ref{cor:intro-nnls}}
\label{app:nnls}

It has been noted before~\cite[Section~1.1]{DHR} that the active-set algorithm of Lawson and Hanson~\cite[Algorithm~(23.10)]{LH74} for nonnegative least squares is \emph{essentially equivalent} to Wolfe's method with the linopt rule.
Unfortunately, we are not aware of a proof or a more precise version of this statement.
Thus, in this section, we give a rigorous proof of the close connection between the two algorithms.
This connection then shows that the lower bound for linopt from the atom of minimum norm (Lemma~\ref{lem:km-start}) carries over to the algorithm of Lawson and Hanson applied to their reduction of least-distance programming~\cite[Chapter~23]{LH74}.

For atoms~$\ivec{p}{j}\in\Q^D$ with~$j\in J$, Lawson and Hanson~\cite[Chapter~23]{LH74} reduce the least-distance problem~$\min\{\norm{\vec{x}}:\vec{x}^\top\ivec{p}{j}\geq1\text{ for all }j\in J\}$ to the NNLS problem
\begin{equation}
    \min\ \norm{E\vec{\nu}-\vec{f}},\quad\text{subject to}\quad\vec{\nu}\geq\zero,
    \label{eq:nnls}
\end{equation}
where~$E\in\Q^{(D+1)\times J}$ has columns~$((\ivec{p}{j})^\top,1)^\top$, for~$j\in J$, and~$\vec{f}\coloneqq(\zero^\top,1)^\top\in\Q^{D+1}$.

Given an NNLS problem of the form \eqref{eq:nnls}, considering an arbitrary matrix~$E$ with columns indexed by~$J$ and an arbitrary vector~$\vec{f}$ for now, the algorithm of Lawson and Hanson operates as follows.
It maintains a set~$F\subseteq J$ and a vector~$\vec{\nu}\geq\zero$ that vanishes outside~$F$, starting with~$F=\emptyset$ and~$\vec{\nu}=\zero$.
An \emph{outer iteration} computes~$\vec{r}=E^\top(\vec{f}-E\vec{\nu})$.
If~$r_j\leq0$ for all~$j\notin F$, the algorithm stops.
Otherwise it adds to~$F$ an index~$j\notin F$ that maximizes~$r_j$, and an inner loop follows.
It computes the vector~$\vec{\zeta}$ that vanishes outside~$F$ and minimizes~$\norm{E\vec{\zeta}-\vec{f}}$.
If~$\zeta_i>0$ for all~$i\in F$, it sets~$\vec{\nu}\gets\vec{\zeta}$, and the outer iteration ends.
Otherwise, in an \emph{inner iteration}, it sets~$\vec{\nu}\gets\vec{\nu}+\vartheta(\vec{\zeta}-\vec{\nu})$ with~\mbox{$\vartheta=\min\{\nu_i/(\nu_i-\zeta_i):i\in F,\ \zeta_i\leq0\}$}, removes from~$F$ every index~$i$ with~\mbox{$\nu_i=0$}, and repeats.
Consider Algorithm~\ref{alg:lh} for a pseudocode in the style of Algorithm~\ref{alg:wolfe}, where~$\vec{r}=E^\top(\vec{f}-E\vec{\nu})$ refers to the current vector~$\vec{\nu}$, and where~$\ivec{\zeta}{F}$ denotes the vector that vanishes outside~$F$ and minimizes~$\norm{E\vec{\zeta}-\vec{f}}$.

\begin{figure}[t]
  \centering\scriptsize
  \renewcommand\mycommfont[1]{\scriptsize\ttfamily\textcolor{gray}{#1}}%
  \setlength{\algomargin}{2pt}%
  \begin{minipage}[t]{.65\linewidth}
    \begin{algorithm}[H]
    \begin{minipage}{\dimexpr\textwidth-13pt}
        \setstretch{1.0}

	    \caption{\algcapstrut\textsc{LawsonHanson}($E$,~$\vec{f}$)}\label{alg:lh}
	    \Input{\algstrut{}matrix~$E$ with columns indexed by~$J$, vector~$\vec{f}$}
	    \Output{\algstrut{}set~$F$, optimal solution~$\vec{\nu}$}
	    \vgap
	    \hrule
	    \vgap
	    \setcounter{AlgoLine}{0}

	    \nl \algstrut$F\gets\emptyset$, $\vec{\nu}\gets\zero$

        \SetInd{4pt}{4pt}
	    \nl \algstrut\While(\tcp*[f]{outer iterations}){$r_j>0$ for some~$j\in J\setminus F$}
	    {
	    \vspace{2pt}
		 \nl \algstrut{}choose~$j\in J\setminus F$ maximizing~$r_j$\tcp*{$\vec{r}=E^\top(\vec{f}-E\vec{\nu})$}
		 \nl \algstrut$F\gets F\cup\{j\}$\tcp*{insertion,~$\nu_j=0$}
		 \nl \algstrut\While(\tcp*[f]{inner iterations}){$\icoord{\zeta}{F}{i}\leq0$ for some~$i\in F$}
		 {
		 \nl \algstrut$\vec{\zeta}\gets\ivec{\zeta}{F}$\tcp*{least squares on~$F$}
		 \nl \algstrut$\vartheta\gets\min\{\nu_i/(\nu_i-\zeta_i):i\in F,\ \zeta_i\leq0\}$

		 \nl \algstrut$\vec{\nu}\gets\vec{\nu}+\vartheta(\vec{\zeta}-\vec{\nu})$\;
		 \nl \algstrut$F\gets F\setminus\{i\in F:\nu_i=0\}$\;
		 }
		 \nl \algstrut$\vec{\nu}\gets\ivec{\zeta}{F}$\tcp*{acceptance}
	    }
	    \nl \algstrut\Return~$F$,~$\vec{\nu}$
	    \vspace{2pt}
    \end{minipage}
	\end{algorithm}
  \end{minipage}
\end{figure}

The key observations that yield the equivalence to Wolfe's method with the linopt rule are the close connections between the vector~$\vec{\zeta}$ and the affine minimizer, and between the entries of the vector~$\vec{r}$ and the gaps, which are given in the following lemma.

\begin{lemma}\label{lem:nnls-ls}
Let~$S$ be an affinely independent support of the atoms~$\ivec{p}{j}\in\Q^D$ with~$j\in J$, and let~$E$ and~$\vec{f}$ be as in~\eqref{eq:nnls}.
Then the columns of~$E_S$ are linearly independent, and the vector that vanishes outside~$S$ and minimizes~$\norm{E\vec{\zeta}-\vec{f}}$ is~$\ivec{\zeta}{S}=\ivec{\bar{\lambda}}{S}/(1+\norm{\vec{w_S}}^2)$.
At~$\vec{\nu}=\ivec{\zeta}{S}$, the vector~$\vec{r}=E^\top(\vec{f}-E\vec{\nu})$ satisfies~$r_j=-g_j(\vec{w_S})/(1+\norm{\vec{w_S}}^2)$ for every~$j\in J$.
\end{lemma}

\begin{proof}
The columns~$((\ivec{p}{i})^\top,1)^\top$ with $i\in S$ are linearly independent because~$S$ is affinely independent.
Let~$\rho\coloneqq\norm{\vec{w_S}}^2$, and let~$\vec{\zeta}$ vanish outside~$S$.
Then, we obtain that~$\norm{E\vec{\zeta}-\vec{f}}^2=\norm{\sum_{i\in S}\zeta_i\ivec{p}{i}}^2+(\one^\top\vec{\zeta}-1)^2$.
If~$\one^\top\vec{\zeta}=0$, this is at least~$1$.
Otherwise, we write~$\vec{\zeta}=s\vec{\lambda}$ with~$s\coloneqq\one^\top\vec{\zeta}$ and~$\one^\top\vec{\lambda}=1$.
As~$\vec{w_S}$ is the affine minimizer of~$S$, we get
\[
    \norm{E\vec{\zeta}-\vec{f}}^2=s^2\Bigl\lVert\sum_{i\in S}\lambda_i\ivec{p}{i}\Bigr\rVert^2+(s-1)^2
    \geq s^2\rho+(s-1)^2\geq\frac{\rho}{1+\rho},
\]
where the last inequality holds since the quadratic function~$s\mapsto s^2\rho+(s-1)^2$ attains its minimum~$\rho/(1+\rho)$ exactly at~$s=1/(1+\rho)$.
Since~$s\neq0$, the first inequality is an equality if and only if~$\sum_{i\in S}\lambda_i\ivec{p}{i}=\vec{w_S}$, and, as the affine coefficients of~$\vec{w_S}$ are unique, if and only if~$\vec{\lambda}=\ivec{\bar{\lambda}}{S}$.
The second inequality is an equality if and only if~$s=1/(1+\rho)$.
So the bound~$\rho/(1+\rho)$ is attained exactly for~$\vec{\zeta}=s\vec{\lambda}=\ivec{\bar{\lambda}}{S}/(1+\rho)=\ivec{\zeta}{S}$.
Since~$\rho/(1+\rho)<1$, the vector~$\ivec{\zeta}{S}$ is the unique minimizer.
At~$\vec{\nu}=\ivec{\zeta}{S}$, we have~$\sum_{i\in S}\nu_i\ivec{p}{i}=\vec{w_S}/(1+\rho)$ and~$1-\one^\top\vec{\nu}=\rho/(1+\rho)$, so
\[
    r_j=1-\one^\top\vec{\nu}-(\ivec{p}{j})^\top\sum_{i\in S}\nu_i\ivec{p}{i}=\frac{\rho-\vec{w_S}^\top\ivec{p}{j}}{1+\rho}=-\frac{g_j(\vec{w_S})}{1+\rho}.\tag*{$\qed$}
\]
\end{proof}

We now obtain the connection between the two algorithms by showing that after every outer iteration, the set~$F$ is the current corral~$S$ and~$\vec{\nu}=\ivec{\zeta}{S}$.

\begin{lemma}\label{lem:lh-wolfe}
Consider the Lawson--Hanson algorithm on~\eqref{eq:nnls} and Wolfe's method with the linopt rule, started at the corral that only contains the index that Lawson--Hanson adds in the first iteration, breaking ties in the same way.
If only one coefficient becomes zero in every minor cycle of Wolfe's method, then both methods add the same indices in the same order, and they remove exactly the same index in each step; so the number of outer iterations of the Lawson--Hanson algorithm that add an index equals the number of corrals of the Wolfe run.
\end{lemma}

\begin{proof}
In the first outer iteration of the Lawson--Hanson algorithm, we have~$\vec{\nu}=\zero$ and~$\vec{r}=E^\top\vec{f}=\one$, so the algorithm can choose the first index~$k$ freely.
Then, by Lemma~\ref{lem:nnls-ls}, the algorithm computes the vector~$\ivec{\zeta}{\{k\}}$, whose only nonzero coordinate has the positive entry~$\icoord{\zeta}{\{k\}}{k}=1/(1+\norm{\ivec{p}{k}}^2)>0$,
so the outer iteration ends with~$F=\{k\}$ and~$\vec{\nu}=\ivec{\zeta}{\{k\}}$, matching Wolfe's start.

Suppose that an outer iteration of the Lawson--Hanson algorithm starts with~$F=S$ and~$\vec{\nu}=\ivec{\zeta}{S}$ for the current corral~$S$ of Wolfe's method.
By Lemma~\ref{lem:nnls-ls},~$r_j>0$ if and only if~$j$ is improving at~$\vec{w_S}$, and the maximizers of~$r_j$ are the minimizers of~$g_j(\vec{w_S})$.
So both methods stop at the same time, and otherwise they add the same index~$j$.

We show by induction over the inner iterations that~$F$ is Wolfe's current support~$S'$ and that~$\vec{\nu}=c\vec{\lambda}$ for Wolfe's current coefficients~$\vec{\lambda}$ and some~$c>0$; initially~$c=1/(1+\norm{\vec{w_S}}^2)$.
The support~$S'$ is affinely independent (Section~\ref{sec:wolfe}), so~$\vec{\zeta}=c'\ivec{\bar{\lambda}}{S'}$ with~$c'>0$ by Lemma~\ref{lem:nnls-ls}.
Thus~$\vec{\zeta}>\zero$ on~$F$ if and only if~$S'$ is a corral, and then both methods accept (Line \ref{line:wf-accept} of Algorithms \ref{alg:wolfe} \& \ref{alg:lh}), with~$\vec{\nu}=\ivec{\zeta}{S'}$.
Otherwise, for~$\vartheta'\in[0,1]$,
\[
    (1-\vartheta')\vec{\nu}+\vartheta'\vec{\zeta}
    =\bigl((1-\vartheta')c+\vartheta'c'\bigr)\Bigl(\vec{\lambda}+t(\vartheta')\bigl(\ivec{\bar{\lambda}}{S'}-\vec{\lambda}\bigr)\Bigr),
    \qquad t(\vartheta')=\frac{\vartheta'c'}{(1-\vartheta')c+\vartheta'c'},
\]
and~$t$ is an increasing bijection of~$[0,1]$ onto itself.
So the two segments (Lines~\ref{line:wf-step} and~\ref{line:wf-move} of Algorithms \ref{alg:wolfe} \& \ref{alg:lh}) have the same signs at corresponding parameters.
In the first inner iteration, we have~$S'=S\cup\{j\}$, the inserted index~$j$ has~$\nu_j=0$, and all other indices of~$F$ have~$\nu_i>0$.
By Lemma~\ref{lem:affine-min}(a), we have~$\vec{w_{S'}}^\top\ivec{p}{i}=\norm{\vec{w_{S'}}}^2$ for all~$i\in S'$ and~$\vec{w_S}^\top\ivec{p}{i}=\norm{\vec{w_S}}^2$ for all~$i\in S$, so
\[
    \norm{\vec{w_{S'}}}^2=\vec{w_{S'}}^\top\vec{w_S}=\sum_{i\in S'}\icoord{\bar{\lambda}}{S'}{i}\,\vec{w_S}^\top\ivec{p}{i}=\norm{\vec{w_S}}^2+\icoord{\bar{\lambda}}{S'}{j}g_j(\vec{w_S}).
\]
Further, since~$g_j(\vec{w_S})<0$, Lemma~\ref{lem:affine-min}(a) gives~$\vec{w_S}\neq\vec{w_{S'}}$, and as~$\vec{w_S}\in\aff(S')$, the uniqueness of the affine minimizer gives~$\norm{\vec{w_{S'}}}<\norm{\vec{w_S}}$; thus we get~$\icoord{\bar{\lambda}}{S'}{j}g_j(\vec{w_S})<0$ and~$\icoord{\bar{\lambda}}{S'}{j}>0$.
In later inner iterations, all indices of~$F$ have~$\nu_i>0$, since the only zero coefficient has been removed.
Hence~$\vartheta$ is the largest~$\vartheta'\in[0,1]$ for which the left side is nonnegative, and~$t(\vartheta)$ is Wolfe's step length.
By assumption, exactly one coefficient becomes zero, so the inner iteration of the Lawson--Hanson algorithm removes the index that the minor cycle of Wolfe's method deletes, and the new vector is again a positive multiple of Wolfe's coefficients.
Induction over the inner and outer iterations completes the proof.
\qed
\end{proof}

The Lawson--Hanson algorithm can choose the entering index in the very first outer iteration freely (first sentence in the proof of the previous lemma).
However, we have only proven a lower bound for linopt when it starts in the atom of minimum norm (Lemma~\ref{lem:km-start}).
Thus, the previous lemma only yields the desired lower bound for the Lawson--Hanson algorithm if it chooses to add the index of the atom of minimum norm in the first iteration.
Instead of assuming a specific tie-breaking rule in the Lawson--Hanson algorithm, in the following proof, we slightly adapt the matrix~$E$ such that the run of the algorithm corresponds to a run on~$E$ that adds the index of the atom of minimum norm first, without a tie.

\begin{proof}[of Corollary~\ref{cor:intro-nnls}]
We show, more precisely, that for every~$n\geq1$, there is a matrix~$E \in \Z^{(n+3)\times(2n+2)}$, of encoding length polynomial in~$n$, on which the Lawson--Hanson algorithm for~$\min\{\norm{E\vec{\nu}-\vec{e_{n+3}}}:\vec{\nu}\geq\zero\}$ performs exactly~$2^n+n+1$ outer iterations that add an index, for every rule for breaking ties.
Take the atoms of Lemma~\ref{lem:km-start} for the lower bound on linopt, and let~$E$ and~$\vec{f}=\vec{e_{n+3}}$ be as in~\eqref{eq:nnls}.
We multiply the column~$\vec{E_{2n+2}}$ by~$2$, and call the result~$\widetilde{E}=E\operatorname{diag}(\vec{\eta})$, where~$\eta_{2n+2}=2$ and~$\eta_j=1$ otherwise.
The substitution~$\vec{\nu}=\operatorname{diag}(\vec{\eta})\widetilde{\vec{\nu}}$ maps the iterates of the algorithm on~$\widetilde{E}$ to iterates on~$E$: it preserves~$E\vec{\nu}=\widetilde{E}\widetilde{\vec{\nu}}$, the vectors~$\vec{\zeta}$ that vanish outside~$F$ and minimize~$\norm{E\vec{\zeta}-\vec{f}}$, the ratios~$\nu_i/(\nu_i-\zeta_i)$, and the sets~$\{i\in F:\nu_i=0\}$ (Line 9 of Algorithm~\ref{alg:lh}), and it multiplies~$r_j$ by~$\eta_j$.
So the algorithm on~$\widetilde{E}$ is the algorithm on~$E$ with the rule that adds an index maximizing~$\eta_jr_j$.
In the first outer iteration, we have~$\vec{r}=\vec{1}$, so~$\eta_jr_j=\eta_j$.
Thus, the algorithm chooses to add the index~$2n+2$, without a tie.
Every corral of the Wolfe run of Lemma~\ref{lem:km-start} contains the index~$2n+2$, so by the proof of Lemma~\ref{lem:lh-wolfe}, the index~$2n+2$ is never again a candidate, and the algorithm on~$\widetilde{E}$ follows the Wolfe run with the linopt rule from~$\{2n+2\}$, which has no ties.
The minor cycles of this Wolfe run are those of Theorem~\ref{thm:equivalence}(b), in each of which only the coefficient of the leaving index becomes zero (proof of Lemma~\ref{lem:pivot}), so Lemma~\ref{lem:lh-wolfe} applies.
By Lemma~\ref{lem:lh-wolfe} and Lemma~\ref{lem:km-start}, the number of outer iterations that add an index is~$2^n+n+1$.
The matrix~$\widetilde{E}$ is integral with~$n+3$ rows and~$2n+2$ columns.
\qed
\end{proof}

\subsection{Proof of Corollary~\ref{cor:intro-pspace}}
\label{app:pspace}

It is PSPACE-complete to predict Dantzig's rule~\cite{FS15}, and this is proved on LPs that are induced by Markov decision processes, so the hardness persists for linopt, even when started at the atom of minimum norm.
Similarly, by~\cite{DS18}, it is NP-hard to determine the number of major cycles of linopt from a suitable corral.

\begin{proof}[of Corollary~\ref{cor:intro-pspace}]
We prove the first part in the following form: given integer atoms with a unique atom of minimum norm and an index~$k$, it is PSPACE-hard to decide whether Wolfe's method with the linopt rule, started at this atom and breaking ties by a fixed order of the indices, ever inserts~$k$, and whether~$k$ is in the final corral.
Fearnley and Savani~\cite[Theorems~1 and~2]{FS15} showed that, for every rule for breaking ties, it is PSPACE-complete to decide whether Dantzig's rule, started at a given feasible basis, ever enters a given variable, and whether a given variable is basic in the optimal basis that it finds.
Their LP~\cite[Appendix~C]{FS15} is the linear program of an MDP (see the proof of Lemma~\ref{lem:rules}) with the right-hand side multiplied by a positive number, which changes no simplex run; their initial basis~$B_0$ is feasible, and Dantzig's rule finds an optimum, so the optimum is finite.
As in the proof of Lemma~\ref{lem:rules}, after scaling the rows to integers and adding the cap, the LP satisfies \ref{ass:integral}--\ref{ass:nondeg} and has the same Dantzig runs from~$B_0$, with the cap added to every basis.
Lemma~\ref{lem:init-run} turns it, in polynomial time, into integer atoms with a unique atom of minimum norm on which Wolfe's method with the linopt rule first inserts the atoms of~$B_0$, the cap, and the shield, and then performs a run of Dantzig's rule, with the same entering and leaving indices; a fixed order of the indices for breaking the ties of linopt induces a rule for breaking the ties of Dantzig's rule.
So an index outside~$B_0$ is inserted if and only if Dantzig's rule enters it, and it is in the final corral if and only if it is in the final basis.

For the second part, Theorem~\ref{thm:equivalence}(c) applies to a single simplex run as soon as all bases of this run are nondegenerate, since Lemmas~\ref{lem:bases-corrals} and~\ref{lem:pivot} are applied only to these bases and the scale~$N$ does not depend on the run.
Disser and Skutella~\cite{DS18} show that it is NP-hard to determine the number of pivots of Dantzig's rule on network-flow LPs whose runs visit only nondegenerate bases and have unique entering indices.
Their instances are minimum-cost flow problems with rational data, solved by the network simplex method, which keeps every arc outside the spanning-tree basis at flow zero or at its capacity~\cite[Section~4]{DS18}; this method is the bounded-variable simplex method, whose bases are the spanning trees and whose reduced costs are those of the arcs~\cite[Section~11.11]{AMO93}.
Their construction contains arcs with infinite capacity~\cite[Section~4.2]{DS18}.
Each of these arcs is incident with the source or sink, which have only outgoing or incoming arcs, respectively.
Let~$F$ be the sum of the positive supplies. 
Flow conservation bounds the flow on each such arc by~$F$. 
We may therefore replace every infinite capacity by~$F+1$: the new upper bounds never become tight, so the network-simplex execution is unchanged (our Lean verification\footref{fn:lean} especially covers this statement). 
These finite capacities have polynomial encoding length.

We bring the instances of~\cite{DS18} into the form~\eqref{eq:lp} by scaling the data to integers, adding a slack variable to every capacity constraint, and deleting one redundant flow-conservation constraint per connected component of the network, so that the constraint matrix has full row rank.
All arc capacities in their instances are positive.
Thus, a spanning tree, together with the slack variables of its arcs, the arcs at capacity, and the slack variables of all other arcs outside the tree, forms a basis of this LP; this basis is nondegenerate if the tree arcs carry flows strictly between zero and their capacities.
The slack variable of an arc at capacity has the negated reduced cost of the arc, so the improving variables, the absolute values of their reduced costs, and the leaving variables agree, and the run of the network simplex method with Dantzig's rule becomes a simplex run with Dantzig's rule on nondegenerate bases.
Adding the cap as in the proof of Lemma~\ref{lem:rules} makes the feasible set bounded without changing this run.
Consequently, it is NP-hard to determine the number of major cycles, that is, of insertions, of Wolfe's method with the linopt rule.
\qed
\end{proof}

\end{document}